\documentclass[12pt]{amsart}

\usepackage[hmargin=2.5cm,bmargin=2.5cm,tmargin=3cm]{geometry}
\usepackage{graphicx}
\usepackage{subfigure}
\usepackage{color}

\usepackage[english]{babel}
\usepackage{amsmath}
\usepackage{amsthm}
\usepackage{mathtools}
\usepackage{amssymb}
\usepackage{algorithmicx}
\usepackage[ruled]{algorithm}
\usepackage{algpseudocode}
\usepackage{parskip}
\algnewcommand{\And}{\textbf{and}}

\theoremstyle{definition}

\newtheorem{thm}{Theorem}[section]
\newtheorem{lem}[thm]{Lemma}

\newtheorem{cor}[thm]{Corollary}

\newtheorem{exmp}{Example}
\newtheorem{remark}[thm]{Remark}

\newcommand{\be}{\begin{equation}}
\newcommand{\ee}{\end{equation}}
\newcommand{\bea}{\begin{eqnarray}}
\newcommand{\eea}{\end{eqnarray}}
\def\var{\mathop{\mathrm{Var}}}

\def\argmin{\mathop{\mathrm{arg\,min}}}
\def\argmax{\mathop{\mathrm{arg\,max}}}

\newcommand{\Graph}{\mathcal{G}}

\newcommand{\cT}{{\mathcal T}}
\newcommand{\cL}{{\mathcal L}}

\newcommand{\cA}{{\mathcal A}}

\newcommand{\cD}{{\mathcal D}}
\newcommand{\cE}{{\mathcal E}}

\newcommand{\cP}{{\mathcal P}}
\newcommand{\cW}{{\mathcal W}}
\newcommand{\cV}{{\mathcal V}}
\newcommand{\cX}{{\mathcal X}}
\newcommand{\cZ}{{\mathcal Z}}

\newcommand\p{\phantom{-}}

\def\var{\mathop{\mathsf{Var}}}

\def\eps{\varepsilon}
\def\Rl{{\mathbb{R}}}

\def\up#1{^{(#1)}}
\let\hat\widehat

\newcommand{\parab}[1]{\paraspace\noindent{\textbf{#1}} }
\newcommand{\paraspace}{\vspace{0.05in}}

\newcommand{\GB}[1]{{\color{blue} (GB: #1)}}

\newcommand{\ND}[1]{{\color{green} (ND: #1)}}

\title[Consistent Merging and Pruning of Subgraphs]{Optimizing
  Consistent Merging
  and Pruning of Subgraphs in Network Tomography}

\thanks{This material is based upon work supported by DARPA and by the National Science Foundation under Grant DMS-1410657 and DMS-1815075. The views, opinions and/or findings expressed are those of the authors and should not be interpreted as representing the official views or policies of the Department of Defense or the U.S. Government.}

\author[M.~Ettehad]{Mahmood Ettehad}
\address{Department of Mathematics, Texas A\&M University, College
  Station, TX 77843-3368, USA}

\author[N.~Duffield]{Nick Duffield}
\address{Department of Electrical and Computer Engineering, Texas A\&M University, College Station, TX 77843, USA}

\author[G.~Berkolaiko]{Gregory Berkolaiko}
\address{Department of Mathematics, Texas A\&M University, College Station, TX 77843-3368, USA}

\begin{document}

\begin{abstract}
A communication network can be modeled as a directed connected graph
with edge weights that characterize performance metrics such as loss
and delay. Network tomography aims to infer these edge weights from
their pathwise versions measured on a set of intersecting paths
between a subset of boundary vertices, and even the underlying graph
when this is not known. Recent work has established conditions under
which the underlying directed graph can be recovered exactly the
pairwise Path Correlation Data, namely, the set of weights of
intersection of each pair of directed paths to and from each
endpoint. Algorithmically, this enables us to consistently fused
tree-based view of the set of network paths to and from each endpoint
to reconstruct the underlying network.

However, in practice the PCD is not consistently determined
by path measurements. Statistical fluctuations give rise to inconsistent inferred weight of edges from measurement based on different endpoints, as do operational constraints on synchronization, and deviations from the underlying packet transmission model. Furthermore, ad hoc solutions to eliminate noise, such as pruning small weight inferred links, are hard to apply in a consistent manner that preserves known end-to-end metric values. 

This paper takes a unified approach to the problem of inconsistent
weight estimation. We formulate two type of inconsistency:
\textsl{intrinsic}, when the weight set is internally inconsistent,
and \textsl{extrinsic}, when they are inconsistent with a set of known end-to-end path metrics. In both cases we map inconsistent weight to consistent PCD within a least-squares framework. We evaluate the performance of this mapping in composition with tree-based inference algorithms.
\end{abstract}

\keywords{Network tomography, end-to-end measurement,
covariance, logical trees, asymmetric routing, unicast probing}


\maketitle

\section{Introduction}\label{sec:intro}
\subsection{Tomography and Graph Fusion}

Network performance tomography seeks to infer edge metrics and even
the underlying network topology by fusing measurements of streams of packets traversing a set of network paths. Abstractly, for additive metrics, a putative solution to the network tomography problem attempts to invert a linear relation between the set of path metrics $\cD$ and the link metrics $\cW$ in the form
\begin{equation}
\label{eq:lin} 
\mathcal{D} = \mathcal{A} \mathcal{W}
\end{equation}
Here $\cA$ is the incidence matrix of links over paths, $\cA_{P,l}$ is 1 if path $P$ traverses link $l$, and zero otherwise. The linear
system (\ref{eq:lin}) is generally underconstrained in real-life
networks, and hence does not admit a unique solution. Mean packet delay and log transmission probabilities are examples of such additive path metrics. In practice, both the measurement and inference functions are distributed across a set of network hosts. Each host performs inference from packet measurement on the subset of network paths of which it is an endpoint (i.e. the source or the destination of measurement packets). The inference produced by each host takes the form of a \textsl{logical} weighted subgraph that estimates the spanning graph of the paths terminating at that host, which each logical edge representing a subpath comprising of one or more edges in the underlying network, with a weight corresponding to the aggregate performance metric on that subpath. Hosts exchange these inferred subgraphs or raw packet measurement data with other hosts, or transmit these to a central location where they are fused to perform network-wide inference.

Therefore, fusion of partial subgraphs is a key task both at
individual hosts and for network wide inference. For example, a common inference primitive involves a host correlating two end-to-end
performance measurements collected from routes to a pair of remote
hosts. The result is a logical spanning binary tree with a single
interior vertex and three leaves (the root host and two other hosts),
see Figure~\ref{fig:tree-SR} (left). The root host then fuses the set of binary trees obtained by iteration over all remote host pairs, in order to infer logical tree spanning paths between itself and the other endpoints \cite{831405}. For network-wide inference, the set of such trees generated by all measurement hosts would be fused to infer the spanning logical performance network that connects
them. 

This program prompts three questions. First: under what conditions is
the network identifiable in the sense that distinct values of network
parameters (topology and edge metrics) can be distinguished in the limit of a large number of packet measurements? Second: what algorithms can identify the network parameters in this limit? Third: how are these algorithms best adapted to work with finite measurement data in the sense of being applicable and performing accurately? The first question has recently been answered for a wide class of inference problems on networks with asymmetric paths between host pairs \cite{10.1088/1361-6420/aae798}. Network level inference is performed by fusing source and destination based trees at each measurement host, characterized by their Path Correlation Data (PCD), namely the weight of the intersection of any two paths that share an origin or destination. Necessary and sufficient conditions for a network graph to be reconstructible from the PCD were established in \cite{10.1088/1361-6420/aae798} together
with an explicit reconstruction algorithm. Thus when the PCD are
identifiable from path pair measurements, the full network is
identifiable under the reconstruction conditions. 

The conditions for reconstruction are: (i) each edge is traversed by at least one path connecting two boundary vertices; (ii) each
non-boundary vertex is \textsl{nontrivial} in the sense that in-degree and out-degree are not both equal $1$; and (iii) each non-boundary vertex $x$ is \textsl{nonseparable} in the sense that the set of paths that pass through $x$ cannot be partitioned into two or more subsets with non intersecting end point sets.

\parab{Practical Inference from Inconsistent Data.}
In distinction to the limiting framework just described, actual
network measurements will not provide a unique value of PCD that is
consistent among distinct path pairs having the same intersection.
Such inconsistencies may have both systematic or statistical
origins and are detailed in Section~\ref{sec:challenge} below. The aim of the present paper is to adapt inconsistent data to the network inference algorithms by constructing empirical PCD with inconsistencies removed according to natural optimality criteria.

\subsection{Network Model and Partial Network Graphs}
We represent the communications network by a a directed edge-weighted
graph $G = (V,E,\cW)$ be with vertices $V$, edge set $E$, a single
edge-based non-negative metric $\cW: E \to \mathbb{R}_{\geq0}$. Edges in $E$ represent links between router identified with vertices, while the weights represent packet performance metrics associated with each edge. We do not assume edge symmetry: $(u,v)\in E$ does not imply $(v,u)\in E$, or weight symmetry: $w_{u,v}$ and $w_{v,u}$ need not be equal. A \emph{partial network graph}
$\Graph=(G,V_B,\cP)$ consists of the graph $G$ together with a set
$V_B$ of boundary vertices and the set $\cP$ of directed paths
$P_{u,v}$ between \emph{some} ordered pairs $(u,v)$ of boundary
vertices. We shall call $u$ the \textsl{source} and $v$ the \textsl{receiver} associated with the path $P_{u,v}$. In the context of network tomography, the boundary nodes $V_B$ act as the sources and sinks of measurement packets that traverse the network on the paths in $\cP$. The remaining vertices $V_I=V\setminus V_B$ will be called the interior vertices. While the paths may be the smallest-weight paths through the network, in general we do not make this assumption. 

Our standing assumptions concerning $\cP$ are
\begin{itemize}
\item[(i)] \textit{uniqueness}: for any given pair $(u,v)$ in $V_B$, there is at most one path $P_{u,v}\in\cP$ connecting them in that direction.
\item[(ii)] \textit{no interior boundaries:} no $v\in V_B$ is an interior node of any path in $\cP$. A partial network without this property can be made conformant by replacing such $v$ by an interior node to which it is connected with a zero weight edge; see e.g. Figure 5 in \cite{10.1088/1361-6420/aae798}).
 \item[(iii)]  \textit{path consistency}: if two vertices $u$ and $v$
appear in two paths $P_{b_1,b_2}$ and $P_{b_1',b_2'}$ in the same order, then the subpaths connecting $u$ and $v$ in $P_{b_1,b_2}$ and
$P_{b_1',b_2'}$ are identical (Note that path consistency implies tree consistency of \cite{10.1088/1361-6420/aae798}).
\end{itemize}
The performance metrics are considered additive in the sense that the performance metric of the path $P_{u,v}$ is the sum $W_{u,v}=\sum_{e\in P_{u,v}}w_e$. The same notation extends to sums over subpaths that terminate at interior vertices. Examples of additive performance metrics include packet delay, and negative log transmission probability. 

The notion of partial network graph will be used to describe both the underlying network (in which case it is often not ``partial'': there is a path in $\cP$ between each ordered pair of boundary vertices) and various inferred networks (in which case it is ``partial'' due to the limited data available to the agent preforming inference). In the latter case, the topology $G = (V,E, \cW)$ may be known in advance, or it may itself be inferred from the measurements.

An important example of a partial network graph is the \emph{source tree} $\cT_{b}^S$ of a boundary vertex $b\in V_B$.  This is the minimal subgraph supporting the path set $\cP(\cT_{b}^S) = \{P_{b,v} \colon v\in V_B\}$. Correspondingly, the \emph{receiver tree} $\cT_{b}^R$ is induced by the path set $\{P_{v,b} \colon v\in V_B\}$.

\subsection{The Challenge of Subgraph Weight Inconsistency}\label{sec:challenge}

A key challenge for fusion derives from the \textsl{weight inconsistencies} between different inferred subgraphs of the
network.  We say that two partial network graph $\Graph\up 1$ and $\Graph\up 2$ have inconsistent weights if there is a pair of vertices $u,v \in V \up 1_B \cap V \up 2_B$ such that $P_{u,v} \in \cP \up1$ and $P_{u,v} \in \cP \up2$, yet the corresponding path weights are unequal: $W\up 1_{u,v}\ne W\up 2_{u,v}$. Causes of inconsistencies include:
\begin{itemize} 
\item [($\text{C}_1$)] \textbf{Statistical variation:} when inference of two partial network graphs yields different weight estimates for edges in their intersection.

\item [($\text{C}_2$)] \textbf{Imperfect temporal alignment:} some inference methods use time series of path performance metrics over a sequence of time slots. However, the slots for different time series may not be synchronized among different hosts due to variable transit time or clock skew. Although packet sequence numbers can be used to coordinate slot alignment at different hosts, these may not be available for measurements of background traffic.

\item [($\text{C}_3$)] \textbf{Deviations from a model:} violation of
  assumptions, such as edge independence, cause inconsistent
  apportioning of path performance amongst the path's constituent
  edges.
\end{itemize}
Two main consequences of weight inconsistency are as follows.

\begin{figure}[h]
	\begin{center}
		\includegraphics[width=0.5\textwidth]{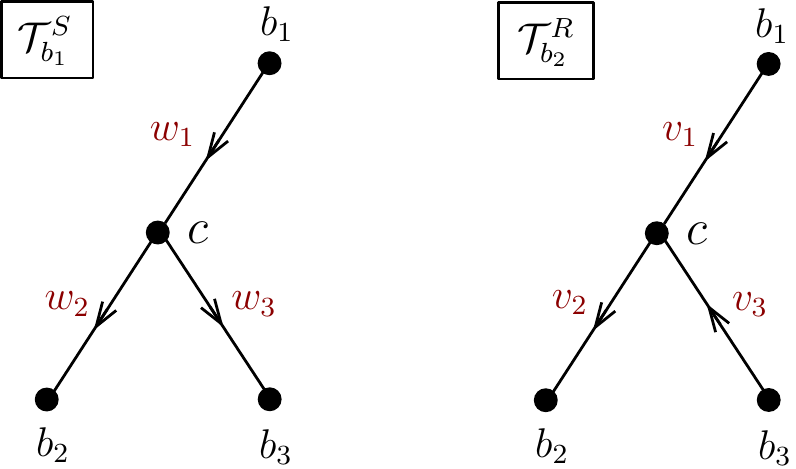}
		\caption{Source tree $\cT_{b_1}^S$ rooted at vertex $b_1$, receiver tree $\cT_{b_2}^R$ rooted at vertex $b_2$ both having version of direct path $P_{b_1,b_2}$.}\label{fig:tree-SR}
	\end{center}
\end{figure}

\parab{Difficulties in Merging Weighted Graphs.} Lack of consistent
path weights prevents application of fusing algorithms that require
them \cite{10.1088/1361-6420/aae798}. To illustrate this, consider in
Figure~\ref{fig:tree-SR} two weighted simple binary trees
$\cT_{b_1}^S$ and $\cT_{b_2}^R$ that have been inferred by primitive
inference using distinct packet time series involving the same three
boundary vertices $b_1$, $b_2$ and $b_3$ with central vertices $c$ and $c'$ respectively. $\cT_{b_1}^S$ is a source tree rooted at vertex $b_1$ while $\cT_{b_2}^R$ is a receiver tree rooted at vertex
$b_2$. The trees have the common directed path $P_{1,2}$ which we wish to merge. However, unless the respective paths weights $w_1+w_2$ and $v_1+v_2$ in $\cT_{b_1}^S$ and $\cT_{b_2}^R$ are equal, it is not immediately clear how to assign interior vertices and edge weights in a merged graph. A naive approach is to assign to the path some combination of the available weights, such as the arithmetic mean. However, changing the paths weights necessitates changing the weights in the constituents edges, which may impact the consistency of other paths that contain those edges, and so on. Instead, what is needed is a principled approach to removing all path inconsistencies by coordinated adjustment of their constituent edge weights.

\begin{figure}[h]
	\begin{center}
		\includegraphics[width=0.5\textwidth]{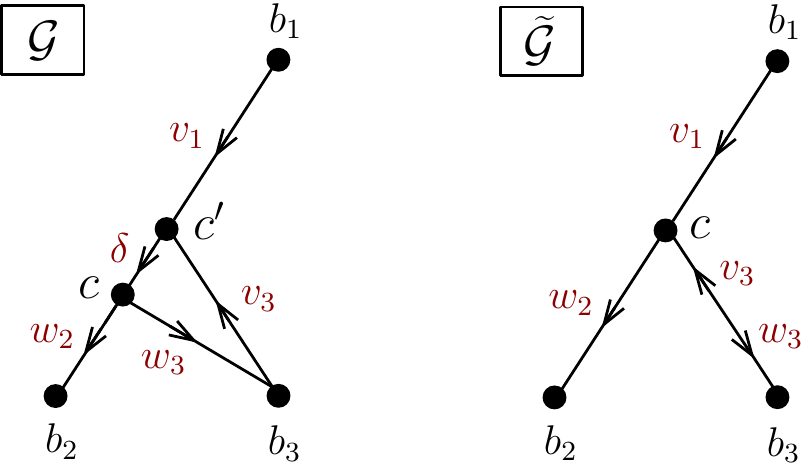}
		\caption{Left: merged graph $\Graph$ under equal path weights $w_1+w_2=v_1+v_2$ for directed path $P_{1,2}$. WNLOG assume $v_1<w_1$ and hence $v_2>w_2$. Right: Pruned graph $\widetilde \Graph$ after removal of edge $(c,c')$ from $\Graph$.}\label{fig:tree-MP}
	\end{center}
\end{figure}

\parab{Topological Noise and Edge Pruning.} Now suppose the weights on path $P_{1,2}$ \textsl{are} consistent, i.e., $w_1+w_2=v_1+v_2$. Then the trees $\cT_{b_1}^S$ and $\cT_{b_2}^R$ may be merged, forming the graph $\Graph$ in Figure~\ref{fig:tree-MP}. Without loss of generality we have assumed $v_1<w_1$ and hence $v_2>w_2$. To merge, the central vertices $c$ from $\cT_{b_1}^S$ and $c'$ from $\cT_{b_2} ^R$ are inserted in the path joining $b_1$ and $b_2$ according to the edge weights in their respective topologies; see
e.g. \cite{10.1088/1361-6420/aae798}. This results in a directed edge $(c,c')$ of weight $\delta=w_1-v_1=v_2-w_2$. The distinction between vertices $c$ and $c'$ may represent asymmetric routing in the underlying network.  However, estimation errors of the type ($\text{C}_1$) will be manifest in the form of extraneous small weight edges.

In order to simplify the inferred topology, such edges are
\textit{pruned} i.e.\ removed from the topology and their endpoints
identified. Criteria for pruning include (a) an edge weight being
less that a threshold performance metric values of interest, and (b)
an edge weight being statistically indistinguishable from zero, e.g.,
on account of being less than some multiple of its estimated standard
deviation. Pruning the edge $(c,c')$ from topology $\Graph$ gives rise to the pruned topology $\widetilde \Graph$ in Figure~\ref{fig:tree-MP}. However, pruning introduces new inconsistencies since the weight $v_1+w_2$ on path $P_{1,2}$ is less than the measured weight $w_1+w_2=v_1+v_2$. Naive approaches to restoring consistency, such as allocating a total weight $w_1+w_2$ in proportion to the weights $v_1$ on edge $(b_1,c)$ and $w_2$ in edge $(c,b_2)$ in $\widetilde \Graph$ will cause weight inconsistencies with other paths containing those edges. Again, we seek a principled way of redistributing the weights of pruned edges while maintaining measured path weights.

\subsection{Problem Statement}
\label{sec:problem_statement}

Motivated by the problems inherent in graph merging and graph pruning described above, we abstract two variant problems in treating inconsistency, as follows:

\parab{Extrinsic Consistency} Let $\Graph=(G,V_B,\cP)$ be a partial network graph. A \textit{target path weight set} $\cZ$ is positive function on a path subset $\widetilde \cP\subseteq \cP$.  We call edge weights $\widetilde\cW$ on $E$ \textit{extrinsically
consistent} with $\cZ$ if $Z_{u,v} = \widetilde{W}_{u,v}$ whenever
$(u,v)\in\widetilde\cP$. Given original edge weights $\cW$ and target path weights $\cZ$, our problem is to find edge weights $\widetilde\cW$ that are (a) extrinsically consistent with the target path weights $\cZ$; (b) close to $\cW$ in a sense to be defined, and (c) readily computable.

\parab{Intrinsic Consistency} Let $\{\Graph\up1,\ldots \Graph\up k\}$
be a set of inferred partial network subgraphs with boundary sets
$V_B\up i \subset V_B = \bigcup_i V_B\up i$, and path sets $\cP\up i$
joining ordered pairs of vertices in $V_B\up i$. Apart from the
boundary points, the graph elements are distinct with no
identification between internal vertices $V\up i\setminus V_B\up{i}$
and edges $E\up i$ from different subgraphs $\Graph\up i$.  We call a
set of weights $\{\widetilde\cW\up 1, \ldots, \widetilde\cW\up k\}$ on the graphs $\{\Graph\up1, \ldots, \Graph\up{k}\}$
\textit{intrinsically consistent} if for any ordered pair
$(u,v)\in V_B$ and all $i$ for which there is a path
$\cP\up i_{u,v} \in \cP\up i$ connecting $u$ and $v$, the path weight
$\widetilde \cW\up i_{u,v}=\sum_{e\in \cP\up i_{u,v}}\widetilde{w}\up
i_e$ is independent of $i$. Given the original set of weights
$\{\cW\up i\}$ on $\{\Graph\up1,\ldots, \Graph\up k\}$, our problem is to find a set of weights $\{\widetilde \cW\up i\}$ that are (a)
intrinsically consistent; (b) close to $\{\cW\up i\}$ in a sense to be defined, and (c) readily computable.

\subsection{Contribution and Outline}

The contributions of the paper are as follows.
\begin{enumerate}
\item We formulate and solve the problems of producing extrinsically or intrinsically consistent weights in a least-squares framework. We seek to minimize the square differences $\Vert \widetilde\cW-\cW\Vert^2$ of the solution and given edge weights, under consistency and positivity conditions for $\widetilde\cW$. In each case, the solutions are expressed in term of the Moore-Penrose inverse of a generalized routing matrix that expresses the linear constraints between path and edge weights. While on the one hand this is a standard approach to constrained optimization, we must make a careful examination of invertibility properties in order to provide a computable solution and avoid singularities. While several
approaches exist in the literature to ensure positivity conditions, we are able to take a simpler approach informed by the network context that small or negative weights are uninteresting and may be set to zero and thereafter ignored.
  
  \item We show how our solution for extrinsic consistency applies to the problem of pruning topological noise.
  
  \item We show how our solution for intrinsic consistency applies to the problem of preparing inconsistent inferred trees for fusion in networks with asymmetric routing.
  
  \item\label{item:comp} We provide a composite application of our methods to the problem of network inference in five stages: (a) merging primitive binary trees at a single root (b) removing noise from the resulting binary trees by pruning (extrinsic consistency); (c) rendering trees with different roots intrinsically consistent; (d) merging trees using the methods of \cite{10.1088/1361-6420/aae798}; (e) removing noise from the merged network by pruning (extrinsic consistency); 
  
  \item We evaluate the composite application in a model-based
  simulation.  We compare its performance in correctly identifying
  subset of network paths that experience performance degradation of a common internal edge. The baseline methods are naive pruning
  and a non-optimal averaging method.
\end{enumerate}

The remainder of the paper is organized as follows.
Section~\ref{sec:related} further motivates our work by describing current approaches for tomographic inference of source and receiver based trees using multicast and unicast probes, and reviews existing approaches to subgraph fusion. Section~\ref{sec:extrinsic} presents the main theoretical results of the paper.  In Section~\ref{sec:consistency_general} we formulate and solve the constrained least-square problem for extrinsic consistency and demonstrate its application to pruning noisy edges. In Section~\ref{sec:intrinsic} we formulate the constrained least-square
problem for intrinsic consistency and solve it for the case of fusing
overlapping inferred trees. Section~\ref{sec:signConstraint} adapts
our methods to take account of positivity constraints on edge metrics. In Section~\ref{sec:expt} we describe and the composite application of intrinsic and extrinsic methods outlined in Section~\ref{item:comp}, metrics for evaluation, and outcomes from model simulation experiments over a range of topologies. After concluding in Section~\ref{sec:conc} we provide the proofs of our results in Section~\ref{sec:proofs} and review the application to our problem of iterative methods for optimization under positivity constraints in the Appendix.


\section{Scenarios for Inferring \& Merging Partial Subgraphs, \\
	and Related Works}\label{sec:related}
In this section we describe scenarios for inference and fusion of \textit{logical network subgraphs} in to which our methods apply, and review prior work on the problem of subgraph fusion.

\subsection{Logical Network Subgraphs in Tomography}
The logical network subgraph associated with a partial network graph $(G,V_B,\cP)$ comprises the partial network graph $(G',V_B,\cP')$ with $G'=(V',E',\cW')$ defined as follows: $V'$ is the union of $V_B$ and the branch points of the path set $\cP$ with $V_I$. $(u,v)\in E'$ if $u,v$ lie on a directed path in $\cP$ without another node in $V'$ between them. $W'_{u,v}$ is the aggregate weight of directed edges in $E$ that connect $u$ to $v$ along the unique path in $\cP$ that joins them. $\cP'_{u,v}\in\cP'$ comprises the edges in $E'$ that are subpaths in $\cP_{u,v}$. In the following scenarios, we assume the full physical network $G=(V,E\up 0,\cW)$ together with a maximal set of possible boundary vertices $V\up 0_B$ and routes $\cP\up 0$ between them. For one or more boundary subsets $V_B\subset V_B\up 0$ and their connecting routes $\cP\subset \cP\up 0$,  tomography is used to estimate the logical network subgraph $(G',V_B,\cP')$ associated with the partial network graph $(G,V_B,\cP)$. Our methods then concern how to consistently fuse the resulting set of logical subgraphs.

\subsection{Multicast Tomography}\label{sec:mult} In multicast tomography, a sequence of multicast probe packets is dispatched from a boundary source along a multicast tree. Successful receipt and transmission latency are recorded for each packet at each receiver. Maximum likelihood estimators for loss \cite{796384}, and discretized delay distribution  \cite{Presti:2002:MIN:611408.611413} are computed for logical edges under independence assumptions for loss and delay on edges. If the topology is not known, each logical source tree is recovered by recursive clustering in which
vertices with the largest common path loss or delay weight are identified as siblings \cite{971737}; see Section~\ref{sec:logicalTreeReconst} for further details. Recovery on the edge weights in a known multicast topology by fusing packet level measurement from trees is proposed in \cite{DBLP:conf/sigmetrics/BuDPT02}. Multicast inference exploits the inherent correlation of multicast packets, each of which occurs once per edge, with copies propagated from each branch point. In order to avoid the requirement for multicast \cite{916283} proposed probing using sets back-to-back unicast packets sent to
pairs of receivers. Their experience approximates that a multicast packet since each unicast packet in a set experiences similar performance on their common path portion.

\subsection{Binary Covariance-Based Tomography}\label{sec:cov}
Let independent random variables  $\cX=\{x_e: e\in E\}$ be associated with each edge and set $X_P=\sum_{e\in P}x_e$ for any subpath $P$. Then
\be\label{eq:cov}
\text{Cov}(X_P,X_{P'})=\text{Var}(X_{P\cap P'})
\ee
and hence any unbiased estimator of the path covariance $\text{Cov}(X_P,X_{P'})$ is an unbiased estimator of the additive path metric $\text{Var}(X_{P''})$ where $P''=P\cap P'$. This approach was first proposed for multicast probing in \cite{832532} where which each packet corresponds to an instance of $\cX$. A different interpretation considers each instance of $\cX$ to represent the edge metric values in force during a time slot. Any flow of packets traversing edge $e$ in a given time slot is assumed to experience performance governed by the same metric value. Hence the window-average performance experienced by distinct flows of packets (even unicast) will be correlated, and more so for larger packet set.  Assuming the instances of $X$ are drawn
i.i.d. over different time slots, then (\ref{eq:cov}) allows estimation of $\var(X_P)$ from measurement of unicast packet on distinct end-to-end paths with intersection $P$; \cite{Nguyen:2007:NLI:1298306.1298339}.

\subsection{Tree Reconstruction from Binary Primitives}
\label{sec:logicalTreeReconst}
As briefly referred to in Section~\ref{sec:mult}, an unknown logical tree can be estimated by recursive clustering on leaf vertices based on largest estimated metric on their common path from the root; the same method can be applied to covariance-based estimates of Section~\ref{sec:cov}. Denote by $X_{k,i}$ the data associated with packet measurements along the path $\cP_{k,i}$ and and estimate $\hat m(X_{k,i},X_{k,j})$ of the aggregate metric the intersection $\cP_{k,i}\cap \cP_{k,j}$. The pair $i,j$ of maximal $\hat m$ are identified as siblings with common parent denoted $\{i,j\}$ while $i$ and $j$ are removed from further consideration. A merged measurement $X_{k,\{i,j\}}$ is then associated--in a metric dependent manner--with the parent, and the process performed recursively until the root is reached. In multicast loss inference $X_{k,i}$ comprises a bitmap indicating which packets reached $i$ from $k$, while for general vertex clusters $A$ and $B$, $X_{k,A\cup B}=X_{k,A}\wedge X_{k,B}$; \cite{796384}. For covariance based estimators $X_{k,A\cup B}$ is a convex combination of $X_{k,A}$ and $X_{k,B}$, with coefficients chosen e.g. to minimize estimation variance \cite{832532}. Due to statistical node, a non-binary node in the underlying network is typically resolved as a set if binary nodes which may then be amalgamated by pruning edges of small inferred weight; see Section~\ref{sec:prune} below.

\subsection{Subgraph Fusion} Prior works have investigated the problem of how to fuse topographically inferred subgraphs. \cite{6731590} fuses quartets (2-source 2-receiver inferred subgraphs) which have been shown to be sufficient to characterize an $M$-source, $N$-receiver network from which they are derived \cite{1354575}. However, information concerning metric edge weights is not exploited in this approach, and so the issue of consistency between different measurements of the same path does not arise. These papers were principally concerned with fusing measurements obtained with striped unicast probes, while our present work is largely agnostic on the measurement mechanism.
The Occam system \cite{2018arXiv180603542S} exploits the idea of binary tree primitives to form source based trees that are then fused by seeking a solution of an optimization problem that seeks a network compatible with the source trees that minimizes over the number of edges and host to host distances. Most recently \cite{10.1088/1361-6420/aae798} proposed fusing source and destination based trees derived from binary tree primitives using passive traffic measurements The key idea is that each path $\cP_{u,v}$ from boundary vertex $u$ to $v$ occurs in both the source tree rooted at $u$ and the receiver tree rooted at $v$. This allows placement of interior vertex on the path according to metric values. However, consistent placement requires equality of the total weight $W_{u,v}$ in each tree, a property that does not hold for estimated weights estimate. Earlier work \cite{Coates:2003:MLT:948205.948230} used detailed timing and packet order information from multiple probe sources to infer overlaps between sourced based trees. Networking tools such as ping and traceroute tools can in principle be used to detect the presence of a given responsive router on distinct paths, although non-responsiveness in encrypted networked and ambiguity due to distinct interface address often limit the utility of this approach.  

\subsection{Tree Pruning}\label{sec:prune}
As discussed above, pruning of low weight edges has been proposed to
reduce topological noise arising from statistical measurement variability \cite{831405}.
If path weight to be preserved, wow should the weight of pruned edges be be ascribed amongst the remaining edges. One approach is to recompute tomographic weights on the pruned topology.
This is well suited to multicast based inference since edge weight estimators exist in non-binary trees \cite{796384} and networks \cite{DBLP:conf/sigmetrics/BuDPT02}.
However, this approach does not generalize to unicast-based network inference. Consider an internal vertex $v$ through which the paths from sources $s_1$ and $s_2$ reach subset $R$ of boundary receivers. The weights of the paths $(s_1,v)$ and $(s_2,v)$ can be individually estimated based on convex combinations of primitive binary estimates from the logical binary tree with vertices $\{\{s_1,v,r_1,r_2\}: r_1\ne r_2 \in R\}$ and $\{\{s_2,v,r_1,r_2\}: r_1\ne r_2 \in R\}$ respectively. However, these estimates will not in general yield equal weights for the edges $\{(v,r): r\in R\}$. This motivates our approach to find weights that are extrinsically consistent (with the path weights) but close to the weights of surviving edges before pruning.


\section{Optimal Assignment of Weights for Extrinsic \\
	     and Intrinsic Consistency}\label{sec:extrinsic}
\newcommand\extpos{\textbf{EXT:POS}}
\newcommand\ext{\textbf{EXT}}
\newcommand\tcW{\widetilde \cW}

\subsection{Extrinsic consistency}
\label{sec:consistency_general}

Let $\Graph=(G,V_B,\cP)$ be a partial network graph and
$\cZ : \cP \to \mathbb{R}_{\geq 0}$ be a target path weight set. If
$\cZ$ is specified on only a subset of $\cP$, we can extend it to all of $\cP$ by assigning $Z_{u,v}=W_{u,v}$ for all other pairs
$(u,v)\in\cP$. We seek a set of weights $\tcW$ on $E$ which is extrinsically consistent with $\cZ$ as a solution to the constrained optimization problem
\begin{align}
  & \min_{\tcW} \Vert \cW - \tcW \Vert^2
    \quad \mbox{such that} \label{eq:extpos:norm}\\
  & \quad \sum_{e\in P_{u,v}} \widetilde{w}_e = Z_{u,v},\
     \forall (u,v) \in \cP \label{eq:extpos:const}\\
  & \quad \widetilde w_e\ge 0,\
    \forall e\in E \label{eq:extpos:pos}
\end{align}
The positivity constraint (\ref{eq:extpos:pos}) originates in the
interpretation that positive edge weights $w_e$ are associated with
performance impairment, while $w_e=0$ indicates no impairment on edge
$e$. We focus first on the optimization problem
(\ref{eq:extpos:norm})-(\ref{eq:extpos:const}) returning to the
positivity constraint in Section~\ref{sec:signConstraint}. The proofs of the following and all other Theorems are given in
Section~\ref{sec:proofs}.

\begin{thm}
  \label{thm:ext}
  Rewrite equation \eqref{eq:extpos:const} as
  \begin{equation}
    \label{eq:matrix_const}
    A\tcW = \cZ,
  \end{equation}
  where $A= \left(A_{(u,v),e}\right)_{{u,v}\in \cP, e\in E}$ denotes
  the incidence matrix of edges in paths,
  \begin{equation}
    \label{eq:A_matr_def}
    A_{(u,v),e} =
    \begin{cases}
      1, &\mbox{if } e \in P_{u,v},\\  
      0, &\mbox{otherwise}.
    \end{cases}
  \end{equation}
  Then the least squares solution of \eqref{eq:extpos:const} which
  minimizes \eqref{eq:extpos:norm} is given by
  \begin{equation}
    \label{eq:ext_sol}
    \tcW = \cW + A'(\cZ-A\cW),
  \end{equation}
  where $A'$ is the Moore-Penrose inverse of $A$.
\end{thm}

\begin{remark}
  We stress that there is no guarantee that \eqref{eq:matrix_const} is consistent (possess at least one solution). If it is not, the
  least squares solution is the best fit with respect to $\ell_2$
  norm.  More precisely, it is the solution $\tcW$ in
  \eqref{eq:ext_sol} is the vector which minimizes
  $\|\cZ - A \tcW\|^2$ and, if there is more than one minimizer,
  $\tcW$ also minimizes $\|\cW-\tcW\|^2$.

  On the practical level, the external constraint is often known only
  approximately. In this situation, it is natural to accept an
  approximate solution to the constraint.
\end{remark}

\begin{remark}
  For numerical computation of the Moore-Penrose inverse, there exist
  simple and stable prescriptions. For example, if $AA^T$ is
  invertible which is equivalent to $Aw=z$ being solvable for any $z$, then $A'$ is given by the formula
  \begin{equation}
    \label{eq:MP_formula}
    A' = A^T (AA^T)^{-1}.
  \end{equation}
  If $AA^T$ is not invertible, formula \eqref{eq:MP_formula} can still be used via regularization or by representing the real symmetric matrix $AA^T$ as $0_K \oplus B$, where $K$ is the null-space of $AA^T$ and $B$ is the restriction of $AA^T$ to the orthogonal complement of $K^\perp$, and interpreting $(AA^T)^{-1}$ as $0_K \oplus B^{-1}$.
\end{remark}

An important special case of Theorem~\ref{thm:ext} is when the
underlying graph is a tree.  In this case \eqref{eq:matrix_const} will have one or more solutions for every $\cZ$. More precisely, let $G$ be a directed tree with root vertex $b$ and leaf set
$\cL_b=V_B\setminus \{b\}$. The following result applies equally to a \emph{source tree}, where there is a unique simple path $P_{b,u}$ for each $u\in \cL_b$ and to a \emph{receiver tree} consisting of unique simple paths $P_{u,b}$, with $u\in \cL_b$.

\begin{cor}\label{thm:inv:tree}
  Let $G$ be a directed tree and let $A$ be the incidence matrix of
  edges over the set of simple paths to (or from) the root vertex
  $\rho$ from (or to) the leaf vertices. Then $A$ has the full row rank
  and the solution that minimizes \eqref{eq:extpos:norm} is given by
  \eqref{eq:ext_sol} with $A'$ computable using \eqref{eq:MP_formula}.
\end{cor}

\subsubsection{An application: edge weight adjustment after pruning}
\label{sec:consistency_pruning}

To show how the general framework of
Section~\ref{sec:consistency_general} applies to pruning we address
the question: how should the weight of pruned edges be assigned to
remaining edges in order to preserve end-to-end path weights.
Consider a partial network graph $G$ and denote by $A$ the incidence matrix of edges over paths, equation~\eqref{eq:A_matr_def}. Pruning an edge amounts to deleting the corresponding column from $A$ and contracting the edge in the underlying graph. All paths in $\cP$ remain connected upon identification of the endpoints of the deleted edge and no further adjustment of $A$ is needed.

Denote by $\tilde E\subset E$ the reduced set of edges and by
$\tilde A$ the corresponding incidence matrix of the pruned
network. Given original edge weights $\cW=\{x_i:i\in E\}$ we seek edge weights $\tcW = \{\tilde w_i: i\in \tilde E\}$ of the pruned network that reproduce the same total weight on all paths in $\cP$,
\begin{equation}
  \label{eq:constraint_after_prune}
  \tilde A \tcW = A \cW,
\end{equation}
and minimize the square distance on the remaining edges,
\begin{equation}
  \label{eq:weights_distance_pruned}
  \Vert \tcW - \cW_{\tilde E} \Vert^2_2
  = \sum_{i\in \tilde E} |\tilde w_i-w_i|^2,
\end{equation}
where $\cW_{\tilde E}$ denotes the restriction of $\cW$ to $\tilde E$.

\begin{cor}\label{cor:prune}
  The least squares solution $\tcW$ of
  \eqref{eq:constraint_after_prune} which minimizes
  \eqref{eq:weights_distance_pruned} is given by
  \begin{equation}
    \label{eq:solution_prune}
    \tcW = \cW_{\tilde E}
    + \tilde A' \left(A\cW -\tilde A\cW_{\tilde E}\right)
  \end{equation}
  where $\tilde A'$ is the Moore-Penrose inverse of $\tilde A$.  If
  the graph $\tilde G$ obtained after pruning is a tree, the matrix
  $\tilde A$ has full row rank.
\end{cor}

\begin{remark}
	The result in Corollary \ref{cor:prune} does not guarantee that all elements of $\tcW$ are positive. In subsection \ref{sec:signConstraint}, we discuss different approaches that can be applied to impose this sign constraint.
\end{remark}

Note that if there is only limited amount of pruning involved, we can
guarantee that system~\eqref{eq:constraint_after_prune} does have at
least one feasible solution even if the graph is not a tree. The
following Lemma sketch the sufficient conditions on the existence of
solution to system~\eqref{eq:constraint_after_prune}.

\begin{lem}
  \label{lem:single_prune}
  Suppose the graph $G$ has no looping edges and a directed edge
  $(u,v)$ has one of the following properties:
  \begin{enumerate}
  \item $u$ is not a boundary vertex and the weight of $(u,v)$ is
    strictly smaller than the weight of any other edges originating
    from $u$, or
  \item $v$ is not a boundary vertex and the weight of $(u,v)$ is
    strictly smaller than the weight of any other edge terminating
    at $v$,
  \end{enumerate}
  then pruning $(u,v)$ results in a consistent system
  \eqref{eq:constraint_after_prune}.
\end{lem}

Since a reasonable pruning scenario is pruning edges of the smallest
weight, one can attempt an iterative application of the above Lemma.
However, there can be problems of topological nature. The following example clarifies the discussion above.   

\begin{figure}[h]
  \begin{center}
    \includegraphics[width=0.7\textwidth]{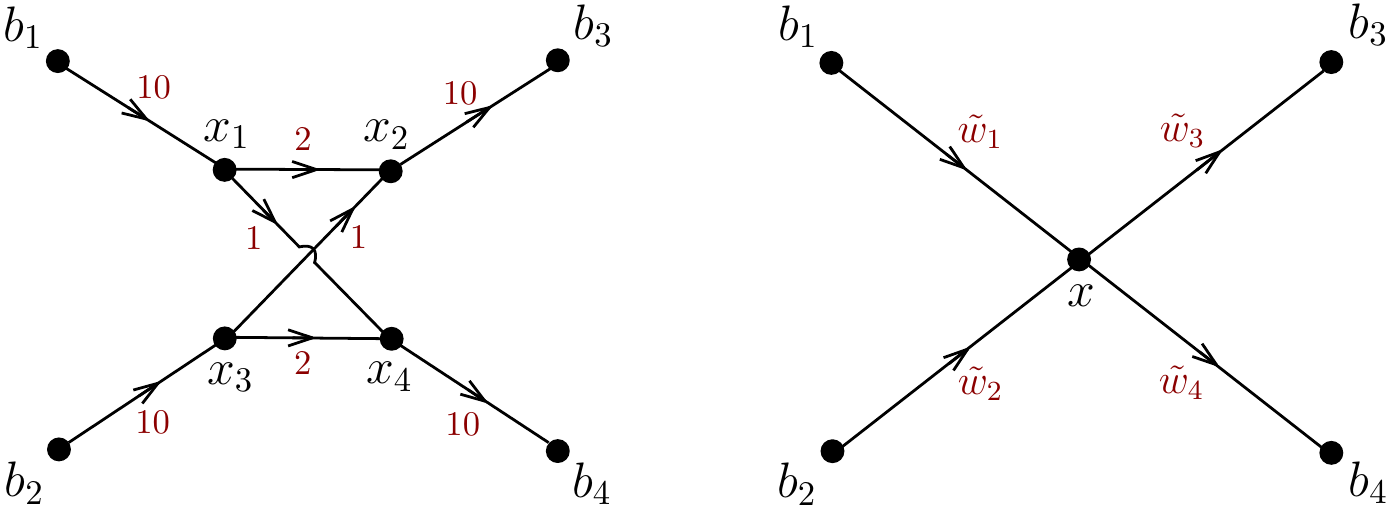}
    \caption{Inconsistency occurs by pruning the edges. Left: the
      original partial graph, right: the pruned
      representation.}
    \label{fig:extrinsicExam}
  \end{center}
\end{figure}
  
\begin{exmp}
  \label{exmp:extInconsistantExample}
  Consider the partial network graph with source $b_1, b_2$ and
  receivers vertices $b_3, b_4$ in Figure
  \ref{fig:extrinsicExam}(left), and the result of pruning four
  shortest edges as shown in Figure \ref{fig:extrinsicExam}(right).
  
  The linear system in equation \eqref{eq:constraint_after_prune} consists of four equations,
  \begin{align*}
    &\tilde w_1 + \tilde w_3 = 22, & &\tilde w_2 + \tilde w_3 = 21,\\
    &\tilde w_1 + \tilde w_4 = 21, & &\tilde w_2 + \tilde w_4 = 22.
  \end{align*}
  Eliminating $\tilde w_1$ from the first column and $\tilde w_2$ from the second column shows the system is inconsistent
  ($\tilde w_3 - \tilde w_4 = \pm 1$).  Applying the result of Corollary \ref{cor:prune} will assign new weights $\tilde w_i = 10.75$ for $i=1,\ldots,4$. While the consistency from end-to-end point of view fails (e.g. $W_{b_1,b_4} \not= \widetilde{W}_{b_1,b_4}$), the
  weights $\tcW$ are optimal in $\ell_2$ sense.
\end{exmp}

\subsection{Intrinsic consistency in Tree-Based Network Inference}
\label{sec:intrinsic}

In the section we address the question of intrinsic consistency (see
Section~\ref{sec:problem_statement}) in the setting of our intended
application, when the graphs $\Graph^{(i)}$ are trees.

More precisely, we assume to be give the set $V_B$ and, for any
$b\in V_B$, two directed trees, a source tree
$\cT^S_b=(V^S_b,E^S_b,\cW^S_b)$ and a receiver tree
$\cT^R_b=(V^R_b,E^R_b,\cW^R_b)$.  In each tree, the root is $b$ and
the leaf set is identified with $V_B\setminus\{b\} =: \cL_b$. All
edges are directed away from $b$ in the source tree $\cT^S_b$ and
towards $b$ in the receiver tree. In applications these trees have
been inferred (with some subsequent pruning) from packet measurements
on an underlying network, with the only known parameter of the network being the set $V_B$.

Note that $u\in \cL_v$ iff $v\in \cL_u$. Each edge $i$ in $E^S_v$
possesses an inferred weight $w^S_{v,i}\geq0$ and likewise each edge
$i$ in $E^R_v$ possesses a inferred weight $w^R_{v,i}\geq0$.  We
emphasize that all trees are distinct with no edge or internal
vertices in common.

For each $u\ne v \in V_B$ let $P^S_{v,u}$ denote the unique path that
connects $v$ to $u$ in $\cT^S_v$.  Similarly, let $P^R_{v,u}$ denote
the unique path connecting $v$ to $u$ in the receiver tree $\cT^R_u$
rooted at $u$.  If the given trees are source and receiver trees
generated from the same underlying network, the paths must be
identical.  In particular, their weights must be the same.

Let $E^*=\bigcup_{v\in V_B}E^R_v\cup E^S_v$ denote the set of all edges over all source and receiver trees and we write
$\cW\in \Rl^{E^*}$ for the vector of their weights, $\cW=\{w^S_{v,i},
v\in V_B,i\in E^S_v\} \cup \{w^R_{v,i}, v\in V_B,i\in E^R_v\}$.
We seek to determine the vector of weights $\tcW =\{\tilde w^S_{v,i}:
v\in V_B,\ i\in E^S_v\}\cup \{\tilde w^R_{v,i}: v\in V_B,\ i\in E^R_v\}\in \Rl^{E^*}$ that minimize the square difference:
\begin{equation}
  \label{eq:sum}
  \Vert \cW - \tcW\Vert^2_2 = \sum_{v\in V_B}
  \bigg(
    \sum_{i\in E^S_v} \left(w^S_{v,i}-\tilde w^S_{v,i}\right)^2
    +
    \sum_{i\in E^R_v} \left(w^R_{v,i}-\tilde w^R_{v,i}\right)^2
  \bigg)
\end{equation}
subject to the common path consistency constraint, 
\begin{equation}
\label{eq:const}
\forall_{v\in V_B},\ \forall_{u\in L_v}:\,
\sum_{i\in P^S_{v,u}}\tilde w^S_{v,i}
= \sum_{i'\in P^R_{v,u}}\tilde w^R_{u,i'}.
\end{equation}

Let $Q$ denote the set of ordered pairs of distinct boundary vertices and let $A$ denote the $|Q|\times |E^*|$ \textsl{signed incidence matrix} of edges over $Q$, defined by
\begin{equation*}
  A_{(v,u),i}=\left\{
    \begin{array}{cl}
      +1 & \mbox{if } i\in P^R_{v,u},\\ 
      -1 & \mbox{if } i \in P^S_{v,u}.
    \end{array}
  \right.  
\end{equation*}
With this notation $(A\cW)_{(v,u)}$ is the asymmetry
between the total weight on the paths from $v$ to $u$ on the receiver
tree with root $u$ and the source tree with root $v$. The constraint
(\ref{eq:const}) is written succinctly as $A\tcW=0$.

\begin{thm}
  \label{thm:intrinsic}
  The matrix $A$ has the full row rank (therefore $AA^T$ is
  invertible) and the solution to the constrained optimization
  (\ref{eq:sum}), (\ref{eq:const}) is given by
  \begin{equation}
    \label{eq:makec} \tcW = \cW - A^T(AA^T)^{-1}A\cW,
  \end{equation}
  with the following error bound
  \begin{equation}
    \label{eq:error_bound}
    \Vert \tcW-\cW\Vert^2_2 \le
    \Vert A\cW \Vert _2^2/2.
  \end{equation}
  This error bound is tight for tree network graphs. 
\end{thm}

\begin{remark}\label{remark:edgePrunHeuristic}
This suggests a possible heuristic for pruning the tree with weights computed using Theorem~\ref{thm:intrinsic}: prune the maximal set of edges of smallest weight, whose sum of square weights does not exceed $\Vert A\cW\Vert_2^2/2$.
\end{remark}

\begin{exmp}\label{exmp:intrinsicExample}
Consider a 3-star network graph with boundary vertices
$V_B=\{a,b,c\}$. For each of the boundary vertices we construct its
source and receiver trees and label by $1$ the edge which appear in two paths, and $2,3$ the edges incidents to boundary vertices in alphabetical order in each of the 6 trees. We order the set $Q$ of pairs of distinct boundary vertices as $\{ab,ac,ba,bc,ca,cb\}$. We order the elements of the six sets of edges $E^R_a,$ $E^R_b,$ $E^R_c,$ $E^S_a,$ $E^S_b,$ $E^S_c$ according to their label. The signed incidence matrix is then given by
\setcounter{MaxMatrixCols}{20}
\begin{equation*}
A =
\begin{pmatrix}
0 & 0 & 0 & 1 & 1 & 0 & 0 & 0 & 0 &  -1 &  -1 & \p0 & \p0 & \p0 & \p0 & \p0 & \p0 & \p0 \\

0 & 0 & 0 & 0 & 0 & 0 & 1 & 1 & 0 & -1 & \p0 & -1 & \p0 &  \p0 & \p0 & \p0 & \p0 & \p0 \\

1 & 1 & 0 & 0 & 0 & 0 & 0 & 0 & 0 & \p0 & \p0 & \p0 & -1 & -1 & \p0 & \p0 & \p0 & \p0 \\

0 & 0 & 0 & 0 & 0 & 0 & 1 & 0 & 1 & \p0 & \p0 & \p0 &  -1 &  \p0& -1 & 0 & \p0 & \p0 \\

1 & 0 & 1 & 0 & 0 & 0 & 0 & 0 & 0 & \p0 & \p0 & \p0 & \p0 & \p0 & \p0 & -1 & -1 & \p0 \\

0 & 0 & 0 & 1 & 0 & 1 & 0 & 0 & 0 & \p0 & \p0 &  \p0 & \p0 & \p0 & \p0 & -1 & \p0 & -1    
\end{pmatrix},
\end{equation*}
which results in
\begin{equation*}
\big(AA^T\big)^{-1}=\frac{1}{90}
\begin{pmatrix}
26&-7&-1&\p2&\p2&-7\\
-7&26&\p2&-7&-1&\p2\\
-1&\p2&26&-7&-7&\p2\\
\p2&-7&-7&26&\p2&-1\\
\p2&-1&-7&\p2&26&-7\\
-7&\p2&\p2&-1&-7&26
\end{pmatrix}.
\end{equation*}

We consider an example in which all edge weights are $1$ except
$w^R_{a,1}$ which is $1+\eps$. Thus $\cW$ is the vector whose first
entry is $1+\eps$ and all other entries $1$. The path weights are
$2$ except for those in the receiver tree rooted at $a$, for which the path weights are $2+\eps$. Then using (\ref{eq:makec}) we compute the consistent edge weights
\begin{equation*}
  \tcW = 1+\frac\eps{90} \big[52,-19,-19,4,-1,5,4,-1,5,2,1,1,14,19,-5,14,19,-5\big]^T.
\end{equation*}
We confirm that, in keeping with Theorem~\ref{thm:intrinsic} above,
$2$ is the smallest eigenvalue of $AA^T$.
\end{exmp}

\subsection{Positivity Constraint Optimization}
\label{sec:signConstraint}
The results in Corollary~\ref{cor:prune} and Theorem \ref{thm:intrinsic} do not guarantee that all elements of consistent weights are positive.  An ad-hoc approach to ensuring positivity would be to prune the edges corresponding to negative $\tilde w$ from the graph(s) and recompute new $\tcW$, iterating until no further negative entries occur in $x$. In numerical experiments a solution satisfying the positivity constraint was achieved in several iterations. There is no guarantee that such naive approach would result in an optimal or nearly optimal solution.

Apart from this naive approach, a more systematic way is to use, for example, path-following method (also known as barrier or interior-point method) for convex quadratic problems which modifies the objective function by adding a nonlinear penalty term with a small coupling coefficient. Let us re-formulate our constrained optimization problem as
\begin{equation}
  \label{eq:primalquad}
  x^* = \argmin_{x \in \mathbb{R}^n} \frac{1}{2} x^T H x + x^T c
  \quad \quad \text{subject to:}
  \quad  Bx = b \quad \text{and} \quad x \geq 0   
\end{equation}
where $H \in \mathbb{R}^{n \times n}$ is positive semidefinite and
$B \in \mathbb{R}^{m \times n}$ has full row rank.
Introducing the \emph{log-barrier term} modulated by a barrier parameter $\mu
\geq 0$ we get 
\begin{equation}
  \label{eq:primalquadProgramBarrier}
  x^* = \argmin_{x \in \mathbb{R}^n} \frac{1}{2} x^T H x + x^T c - \mu \sum_{i=1}^{n} \ln (x_i) \quad \quad \text{subject to:} \quad  Bx = b 
\end{equation}
The idea of the algorithm is that solutions of
\eqref{eq:primalquadProgramBarrier} converge to solutions of
\eqref{eq:primalquad} as $\mu\to0$.  It has been shown that for an
appropriately chosen starting point,
$\mathcal{O}(\sqrt{n} \log\frac{n}{\varepsilon})$ iterations are
required to be $\varepsilon$-close to the optimal solution
\cite{CWZ13}.

Returning to the our problem of imposing the positivity constrains in
\eqref{eq:constraint_after_prune}--\eqref{eq:weights_distance_pruned}
and \eqref{eq:sum}--\eqref{eq:const}, we need to consider two cases:
the matrix $A$ has full row rank or not.

In the context of pruning a tree (Corollary~\ref{thm:inv:tree} and
Corollary~\ref{cor:prune}) or ensuring intrinsic consistency in a set of trees (Theorem~\ref{thm:intrinsic}), we are guaranteed that $A$ has full row rank. In the former case, the standard quadratic
optimization problem can be used by setting $H := 2\mathbb{I}$, $c :=
-2\cW_{\tilde E}$, $B:= \tilde A$ and $b:= A \cW$. In the latter case, the standard quadratic optimization problem can be used by setting $H := 2\mathbb{I}$, $c := -2\cW$, $B:= A$ and $b:= 0$.

With the view of pruning a full reconstructed graph which may not be
a tree (Theorem \ref{thm:ext}), we explain modifications necessary for a matrix $A \in \mathbb{R}^{m \times n}$ which is not of full row rank. In order to circumvent this rank-deficiency, in general, the matrix $A$ can be reduced to full rank matrix via QR factorization or Gaussian elimination with column pivoting \cite{Wright97}. Applying QR factorization applied to the consistent system $A\cW = \cZ$ obtains an $m \times m$ orthogonal matrix $Q$ such that 
\begin{equation}
QA=
\begin{pmatrix}
\bar A \\
0
\end{pmatrix}, \quad 
Q\cZ=
\begin{pmatrix}
\bar \cZ \\
0
\end{pmatrix}
\end{equation}
where $\bar A$ and $\bar \cZ$ have the same number of rows and $\bar A$ has full row rank. Through this construction of matrix $\bar A$, the systems $A \cW = \cZ$ and $\bar A \cW = \bar \cZ$ are equivalent; that is, any vector $\cW$ that satisfies one of these equations also satisfies the other. For detailed example on dealing with rank deficient matrix $A$ and application of QR factorization, see \cite{Wright97}. 


\section{Composite Applications \& Evaluation}\label{sec:expt}
In this section we evaluate the performance of our approach on composite applications necessitating both internal and external consistency.  We generate a set of network graphs (see Section~\ref{sec:setup}), and in each graph simulate the end-to-end packet measurements (details in Section~\ref{sec:data}). The next step is to use the simulated data to infer the logical source and receiver trees rooted at every boundary vertex $b \in V_B$, by applying the recursive clustering approach reviewed in
Section~\ref{sec:logicalTreeReconst}. Low weight edges in the
resulting binary local logical trees are pruned producing non-binary
internal vertices.  In order to keep the set of weights extrinsically
consistent with the end-to-end measurements, the result of
Corollary~\ref{cor:prune} is applied on each pruned tree separately.
The entire set of trees is then made intrinsically consistent using
the results of Section~\ref{sec:intrinsic}. This internal consistency is required for fusion of inferred local trees following the graph reconstruction algorithm of \cite{10.1088/1361-6420/aae798}. The last step before evaluating the inferred graph, is to prune edges with small weight due to topological noises raises from the nature of non-exact measurements. This will be achieved by applying the result of Corollary~\ref{cor:prune} along with positivity constraint.

\subsection{Generating Network Graph}\label{sec:setup}
We first describe a framework on constructing set of random network graphs $\Graph = \{\Graph \up 1, \ldots, \Graph \up N \}$ with each $\Graph\up k = (G \up k, V_B \up k, \cP \up k)$ has desirable characteristics.      

\parab{Random Graphs}
In the first part, we construct random graphs along with their corresponding set of paths among the boundary vertex set through the following steps. For each $k \in [N]$, 
  
\begin{enumerate}
	\item In the first step, a random graph which will be called \textit{underlying graph} $\hat G = (\hat V, \hat E, \hat \cW)$ is constructed with $|\hat V| = m$ arbitrary and each vertex $u \in \hat V$ has degree $\text{deg}(u) = d$ connecting to its randomly selected neighbors. The set $\hat \cW$ is generally asymmetric with respect to edge direction in set $\hat E$. 
	
	\item Among set $\hat V$, $n$ vertices will be selected randomly which represent the boundary set $V_B = \{b_1, \ldots, b_{n}\}$. 
	
	\item The routing between selected boundary vertices $b_i$ and
          $b_j$ i.e. $P_{b_i,b_j}$ with $i \not= j$ follows the
          shortest path with respect to the weights $\hat\cW$. This then form set $\cP$. 
\end{enumerate}
Note that the assigned weight set $\hat \cW$ above is only applied to generate the set of paths $\cP$, and the link's performance characteristic for end-to-end measurements will be assigned in the data generation step below (see subsection \ref{sec:data}). 
\begin{figure}[h]
	\begin{center}
		\includegraphics[width=0.9\textwidth]{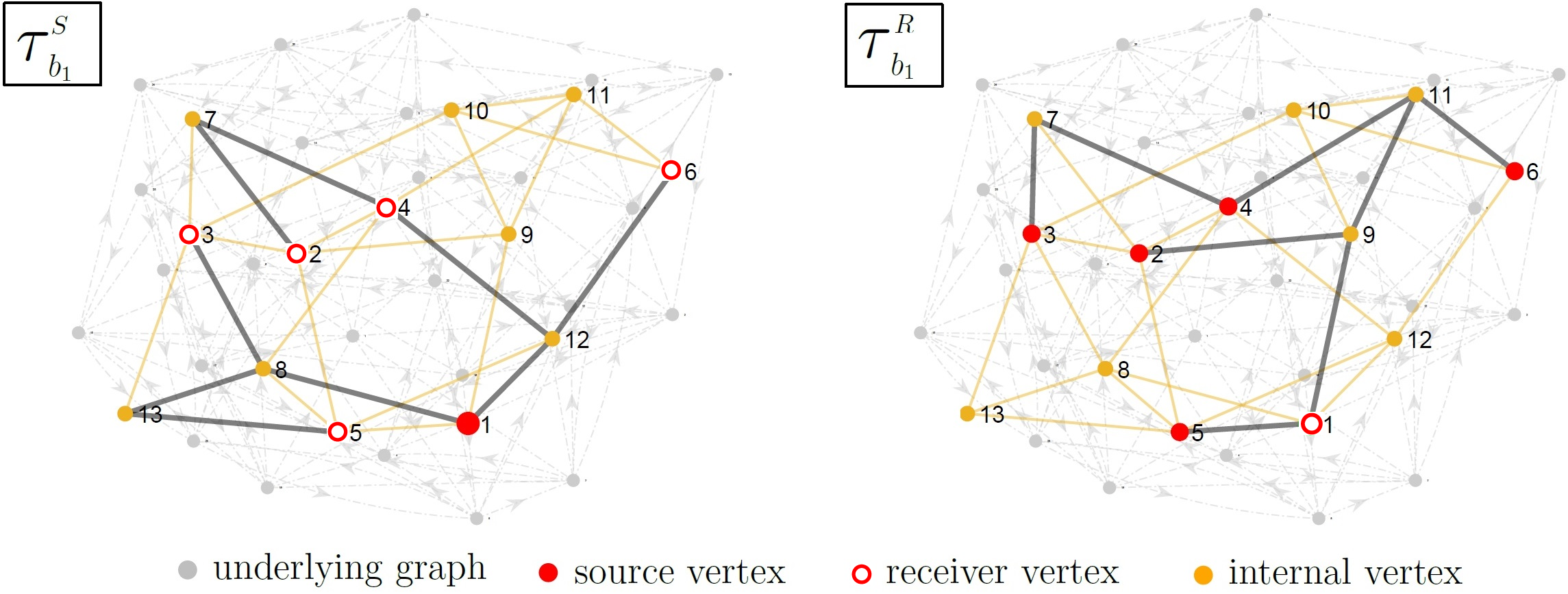}
		\caption{Example on randomly generated source
			$\cT_{b_1}^S$ and receiver $\cT_{b_1}^R$ trees at
			vertex $b_1$. See Section \ref{sec:intro} for the definition of source and receiver trees at vertex $b \in V_B$.}\label{fig:srTrees}
	\end{center}
\end{figure}

An example of application of the three steps discussed above is shown in Figure \ref{fig:srTrees} for graph with six boundary vertices (only the source and receiver trees at vertex $b_1$ are plotted). It should be noted that in the construction above, if the assigned weight set $\hat \cW$ is chosen to be symmetric on edge set $\hat E$, then the communication paths in constructed network graph $\Graph$ will be symmetric, otherwise asymmetric routing will be obtained. 

\subsection{Generating Measurements}\label{sec:data}
\begin{figure}[h]
	\begin{center}
		\includegraphics[width=0.6\textwidth]{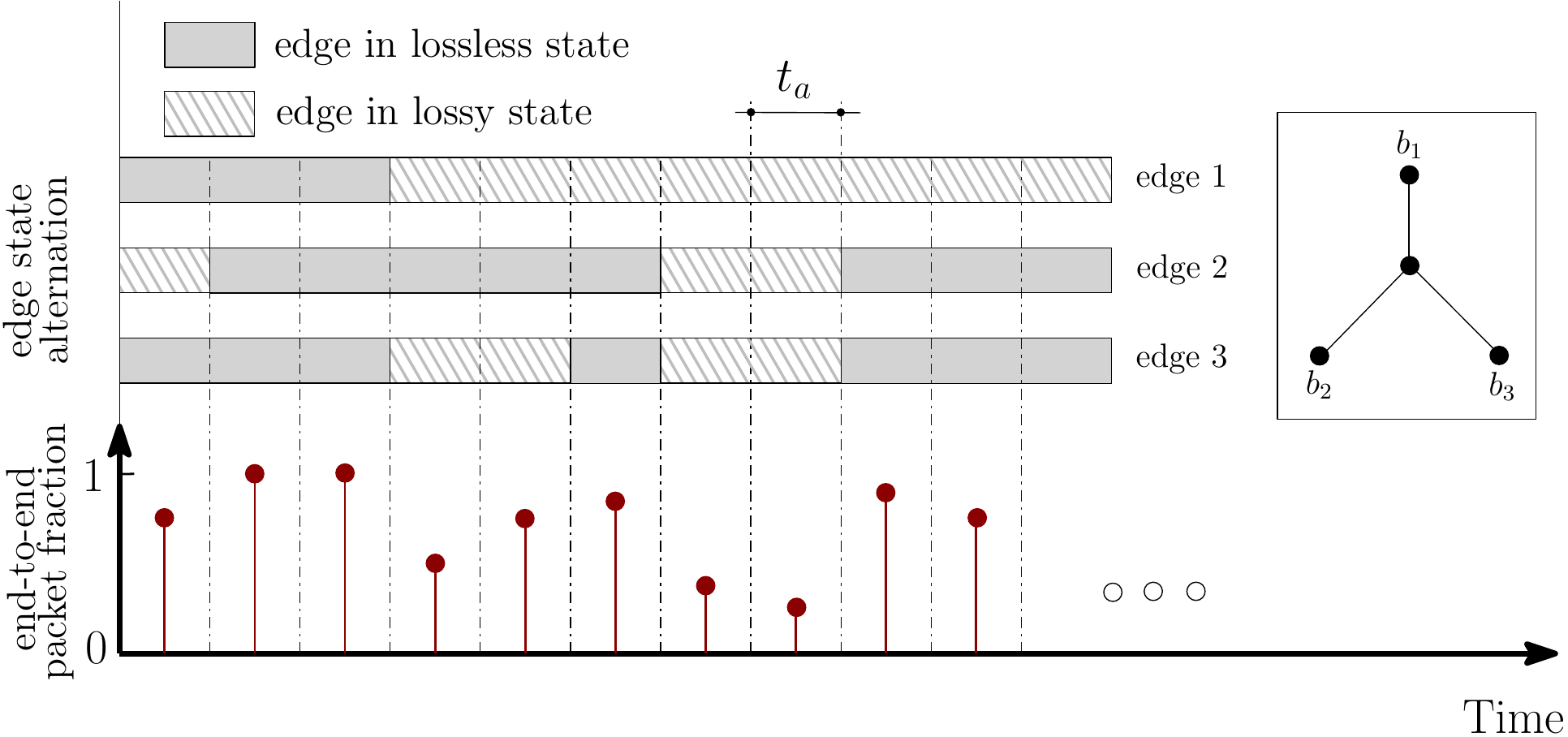}
		\caption{Alternation of edges state and measured
                  end-to-end averaged packet fraction over
                  time.}\label{fig:lossModel}
	\end{center}
\end{figure}

In this part we will discuss the process of generating end-to-end measurements corresponding to set $\cP$ constructed above for network graph $\Graph$. Although there are various ways to generate this set of measurements, in this paper we will follow model-based unicast data acquisition framework described as follows. For graph $\Graph = (G, V_B,\cP)$, each directed edge $i \in E$ will experience alternative lossless and lossy states so that a packet transversing the edge will be dropped with probability zero and $p_l \up i \in \cL$ respectively, where $\cL$ be a set determines the possible edge loss probability. For the experiment with total time $T$, fraction of packets received successfully for selected source and receiver vertex will be recorded over successive averaging windows of length $t_a$ i.e. intervals $((k-1)t_a, kt_a)_{k=1}^N$ with $N t_a = T$.  The time that each directed edge $i$ being in successive lossless or lossy states will be equal to $k_s \up i t_a$ and $k_{\ell} \up i t_a$ with integers $k_s \up i$ and $k_l \up i$ drawn from Poisson distributions with per-determined parameters $\gamma_s$ and $\gamma_l$ respectively. The end-to-end measurements for $P_{b_i,b_j}$ with $b_i, b_j \in V_B$ and $b_i \not = b_j$ then will be time series corresponding to variation of the fraction of successfully received packets over time. The schematic of the above process is shown in Figure~\ref{fig:lossModel} for the network graph with three boundary vertices. 

Finally, the edge weight set $\cW$ for the network graph $\Graph$ is constructed by averaging over the edges performance (e.g. empirical variance or mean of packets loss rate) through the generating  measurements process overt time $T$. In other words, for an edge $i \in E$, the averaging is over all paths $P \in \cP$ so that $i \in P$.    

\subsection{Performance Metric}\label{sec:metrics}
We now detail metric we use to evaluate the closeness of the inferred
network $\hat \Graph$ to the true network $\Graph$. Many metrics have been developed for error tolerant graph matching, including graph edit distance \cite{4767144} and maximal subgraph matching
\cite{42855}. Closer to the present work, \cite{971737} used a notion of receiver matching in tree, namely, that an internal vertex in the true spanning tree associated with a boundary point is well estimates if there is a vertex in the corresponding inferred tree with the same
receiver set.

The motivation for our work is to attribute common origins of
performance degradation on paths, not limited to paths sharing an
endpoint. To this end, we specify a network distance metric
that captures the difference in the extent to which individual edges
influence the performance of boundary-to-boundary transmission.

\parab{Weighted Path Intersection Metric} Given a network graph
$\Graph=(G,V_B,\cP)$, for each directed edge $e\in E$, let $\cV(e)$
denote the set of ordered pairs of boundary vertices
$(u,v)\in V_B\times V_B$ such that $e$ is an edge in the path
$P_{u,v}$ from $u$ to $v$.  We refer to edges $e$ and $e'$ as
equivalent in $\cV(e) = \cV(e')$.  Given another network graph
$\Graph'=(G',V_B,\cP')$ with the identical boundary set $V_B$, define
\begin{equation}
  \label{eq:similar_edges_def}
  \cE(\cV(e),\Graph') := \bigcap_{(u,v)\in \cV(e)} P'_{u,v} \setminus
  \bigcup_{(u,v)\not\in \cV(e)} P'_{u,v}.
\end{equation}
We stress that the set $\cV(e)$ is determined with respect to the
graph $\Graph$ whereas the paths $P'_{u,v}$ on the right-hand side are taken from the graph $\Graph'$.  Thus the meaning of
$\cE(\cV(e),\Graph')$ is the set of edges in $\Graph'$ that are
equivalent (from the point of view of boundary-to-boundary
transmission) to the edge $e$ of $\Graph$. Specifically, they lie in the intersection of paths from $\cP'$ between endpoint pairs $\cV(e)$, but in no other paths in $\cP'$. 
In general, the set $\cE(\cV(e),\Graph')$ may be empty, but if $\Graph' = \Graph$, the set $\cE(\cV(e),\Graph)$ is simply the set of edges equivalent to $e$ as per definition above. 

Let $W$ be a weight function from the set of subsets of $E(G)$ to
non-negative integers; let $W'$ be a similar function on the graph $\Graph'$.
We define
\begin{equation}
  \label{eq:q:wpim}
  Q(\Graph,\Graph')=\frac{ \sum_{e \in E}
    \left|W_{\cE(\cV(e),\Graph)} - W'_{\cE(\cV(e),\Graph')}\right|
  }{
    \sum_{e \in E} W_{\cE(\cV(e),\Graph)}
  },
\end{equation}
where in both summation only one $e$ from each equivalence class is
used. 

The relative error in estimating $W_{\cE(\cV(e),\Graph)}$ is $\left|W_{\cE(\cV(e),\Graph)} - W'_{\cE(\cV(e),\Graph')}\right|/W_{\cE(\cV(e),\Graph)}$.
Thus $Q$ is the average over equivalence classes of edges $e\in E$ of the relative error in estimating $W_{\cE(\cV(e),\Graph)}$, weighted by $W_{\cE(\cV(e),\Graph)}$ i.e., $Q$ is more sensitive to errors in estimating large weights $W_{\cE(\cV(e),\Graph)}$ than for small ones.

The metric assumes that the weights in $\Graph$ and $\Graph'$ to be comparable, e.g., because $W'$ estimates $W$. This is not the case for covariance based estimation as described in Section~\ref{sec:related}. One approach to this is to replace $W$ with weights comparable to $W'$, e.g. by computing edge estimator covariances from a model of packet performance. We here propose two approaches that can be used directly with any weights $W$ and $W'$.
\begin{itemize}
\item [($\text{TM}_1$)] Set $W'_\emptyset=0$ but take
  $W'_{\cE}=W_{\cE} = \sum_{e\in\cE} w_e$ otherwise. In other words,
  we penalize $e$ by its weight if there no edge $e'$ in $\Graph'$
  functionality equivalent to it.
\item [($\text{TM}_2$)] As in $\text{TM}_1$ but with $W_{\cE}=1$ for all
  nonempty subpaths $P$.
\end{itemize} 

\subsection{Numerical Results}\label{sec:numericalResults}
In follows, performance of the reconstruction framework will be evaluated on set of randomly generated network graphs discussed above. The evaluation metric will be based on $(\text{TM}_1)$ and $(\text{TM}_2)$ while evaluation by incorporating edge's weights of reconstructed network will be reported in future work (see the discussion in Section \ref{sec:conc}). In all numerical simulations, we chose the parameters of Poisson distribution which determine the number of averaging windows an edge experiences in lossy or lossless states to be $\gamma_l = \gamma_s = 10$ (the initial state of each edge is chosen randomly). Moreover the loss probability of edge $i$ in lossy state is chosen randomly from  $p_l\up i \in \cL = \{0.05, 0.10\}$ if edge $i$ be in it's lossy state. The number of packets sent form source $b_i \in V_B$ to the set of receiver boundary vertices $V_B \setminus \{b_i\}$ in each averaging window $t_a$ will be fixed and equal to $10^3$ in all simulations. Thereby, number of averaging windows will control the total number of packets sent in unicast mode.  

\begin{figure}[htbp]
	\center
	\subfigure{
		\includegraphics[width=0.45\textwidth]{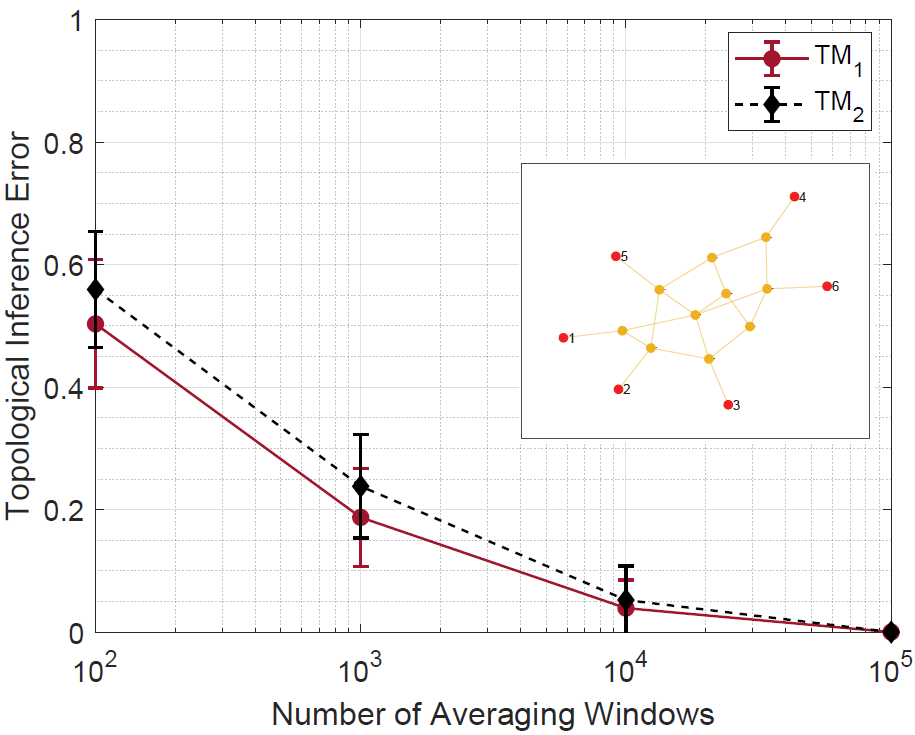}
	}
	\hspace{3mm}
	\subfigure{
		\includegraphics[width=0.45\textwidth]{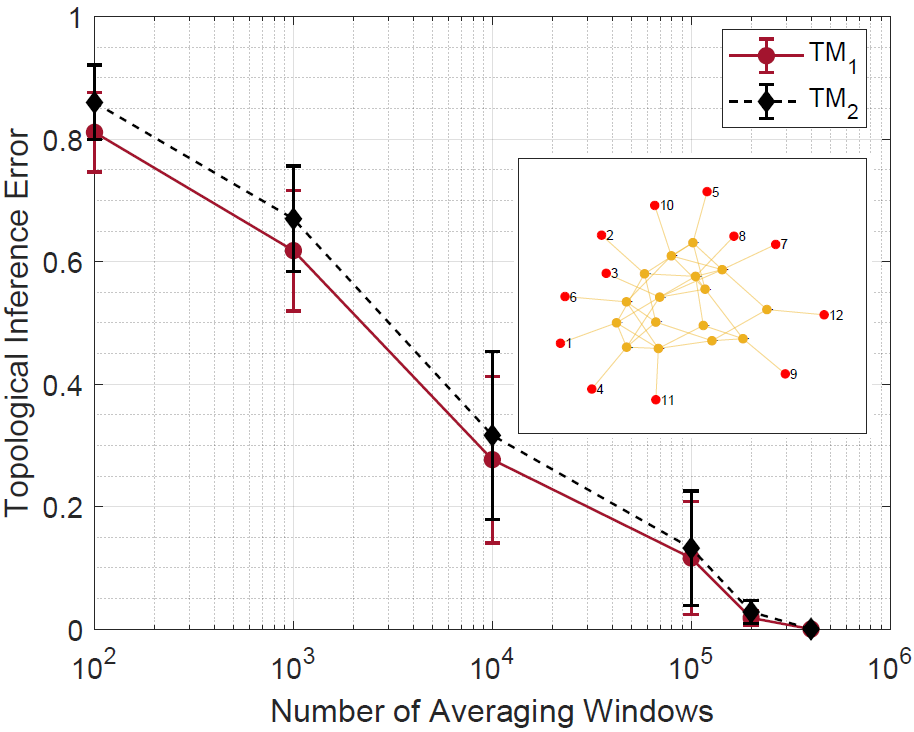}
	}
	\caption{Dependence of topological inference error on the number of averaging windows evaluated based on metric $\text{TM}_1$ and $\text{TM}_2$. Left: network graphs graph with $|V_B| = 6$\ and  right: $|V_B| = 12$.}
	\label{f_ResultsNumWindows}
\end{figure}

Figure \ref{f_ResultsNumWindows} shows the error plot of topological inference over the number of averaging windows for two selected samples of network graph with $|V_B| = 6$ and $|V_B| = 12$ (for each number of averaging windows 50 experiments are conducted). The convergence to full inference accuracy over the number of averaging windows is clear. For the larger network graph (right), higher number of averaging windows is required to achieve same level of inference accuracy happens for six boundary vertex graph (left). The reason is that as the size of network graph increases, the variance of estimators in tree level identification increases. Hence mistaken pairing of non-sibling vertices, or erroneous inclusion or exclusion of vertices in a group is more likely to occur. For same experiment over the network, topological metric $\text{TM}_1$ reports slightly better algorithmic performance compared to $\text{TM}_2$ metric and the difference of these two metrics reduce (converges to zero) over the inference accuracy. 

\begin{figure}[htbp]
	\center
	\subfigure{
		\includegraphics[width=0.45\textwidth]{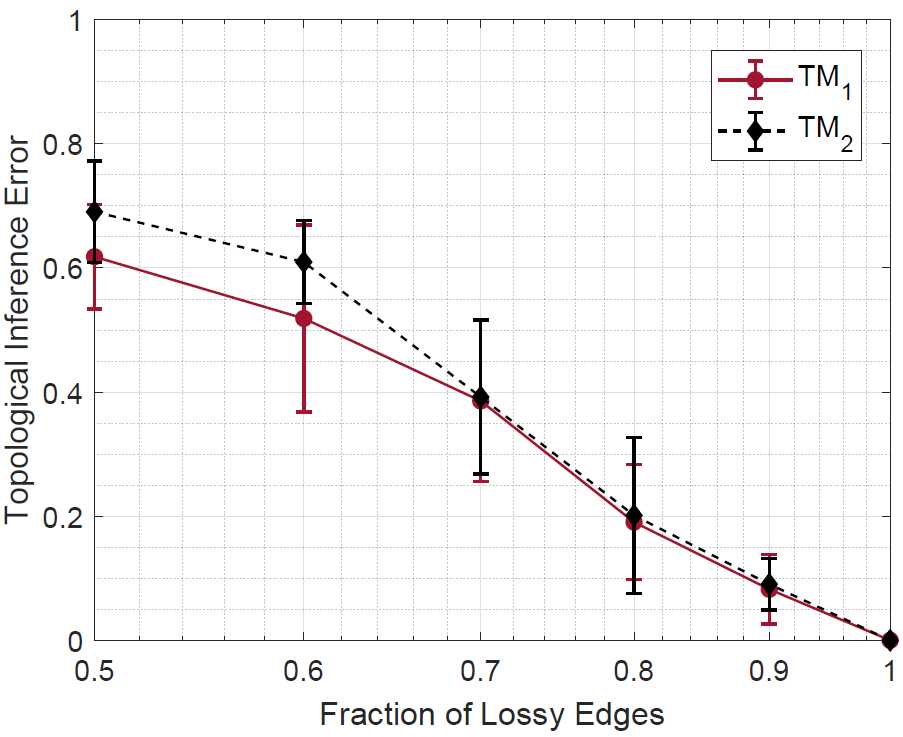}
	}
	\hspace{2.93mm}
	\subfigure{
		\includegraphics[width=0.458\textwidth]{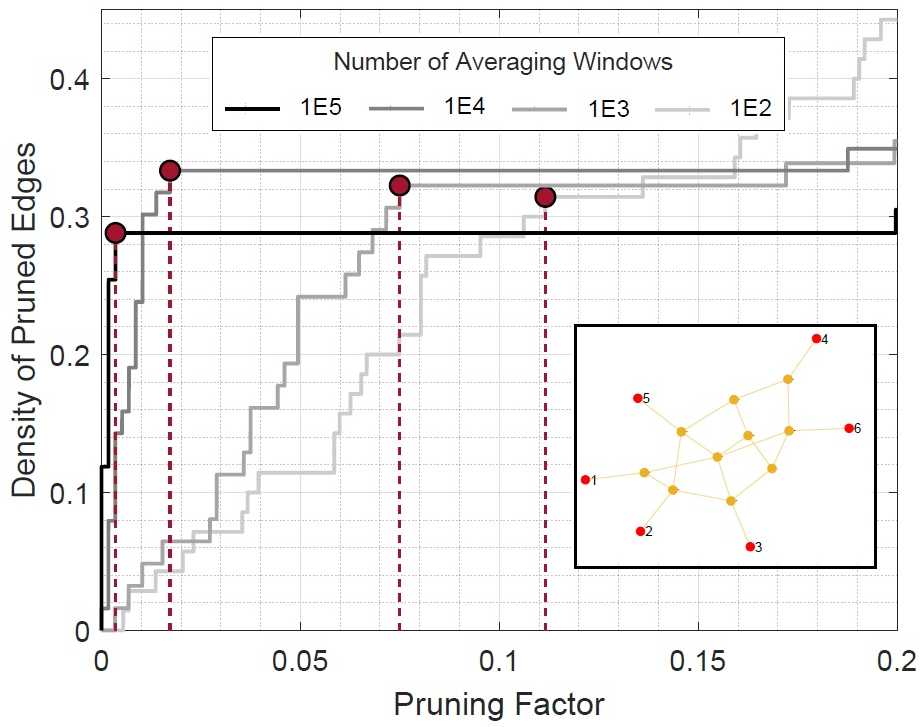}
	}
	\caption{Left: Dependence of topological inference error on
		fraction of lossy edges. Right: Sample on dependence of
		density of pruned edges on pruning factor for different
		number of averaging windows. The vertical dashed lines are the selected pruning factor for the corresponding experiment.}
	\label{f_ResultsFracLoss}
\end{figure}

Figure \ref{f_ResultsFracLoss}(left) shows the error plot of
topological inference over the fraction of edges experiencing lossy
state for 20 randomly selected networks graphs with six boundary vertices. For each network graph and selected fraction of lossy edges, 50 experiments are conducted where edges which experience both lossless and lossy states are selected uniformly in random. Obviously, those edges which are not selected as lossy ones stay in lossless state over the whole experiment. The plot shows that once the fraction of lossy edges increases, the average of inference accuracy approaches to the selected lossy fraction while uncertainty in inference of identifiable subgraph (part of graph constructed by lossy edges) increases for low fraction of lossy edges. The difference of topological inference metrics $\text{TM}_1$ and $\text{TM}_2$ decreases for larger fraction of lossy edges while the two metrics start deviating from each other when the fraction of lossy edges is less than 0.7. 

In the numerical simulations, an edge $i \in E'$ belonging to the reconstructed network graph $\Graph'$ will be pruned if $w'_i/w'_{\max} \leq \delta$ with $\delta \in [0,\delta_{\max}]$ is called pruning factor with $\delta_{\max} \leq 1$. In Figure \ref{f_ResultsFracLoss}(right), the dependence of density of pruned edges on pruning factor ($\delta_{\max} = 0.2$) for different values of number of averaging windows (one sample is selected) is plotted. The result shows that once the number of averaging windows increases then small value of pruning factor is sufficient to prune relatively large number of edges with small weights due to topological noise (see the discussion in Section \ref{sec:extrinsic}) while achieving the same number of pruned edges requires large value of pruning factor for inference based on low number of averaging windows. This is in line with the fact that in the limiting case when the edge weights approaches to the exact value, then the pruning factor required to pruned edges exist due to topological noise approaches to zero. Letting $(\delta_k)_{k=1}^m$ with $\delta_m \leq \delta_{\max}$ to be increasingly ordered point of discontinuities for the plot of dependence of density of pruned edges over pruning factor (e.g. see Figure \ref{f_ResultsFracLoss}(right) for selected number of averaging window), then an ad-hoc approach to choose the pruning factor $\delta^* = \delta_{k^*}$ in each experiment is to find $k^* = \argmax_{k \in [m-1]}(\delta_{k+1} - \delta_k)$. Another possible heuristic approach to prune the inferred trees is to follow Remark \ref{remark:edgePrunHeuristic}.  


\section{Conclusion and Future Work}\label{sec:conc}
In this paper, we extended the method of graph reconstruction to Path Correlation Data in which sampling and measurement noise leads to inconsistencies in path weights reported for different paths. The two main types of inconsistencies due to the nature of non-exact measurements are formulated in the form of extrinsic and intrinsic types, and solved in a least-squares framework. In each case, the solutions are expressed in term of the Moore-Penrose inverse of a generalized routing matrix which express the linear constraints between path and edge weights. This extension provided a unified framework on consistency merging and pruning subgraphs which reduce the problem in network inference through the following five steps: (a) merging primitive binary trees at a single root (b) removing topological noise from the resulting binary trees by pruning (extrinsic consistency); (c) rendering trees with different roots intrinsically consistent; (d) merging trees using the methods of \cite{10.1088/1361-6420/aae798}; and (e) removing topological noise from the merged network by pruning (extrinsic consistency). In this work we discussed the application of prescribed five composite steps on inference of randomly generated network graphs using model based data. Although the proposed metric on inference accuracy is general in the sense that it can capture both fully topological base and weighted based accuracy, but in this work we focused on its topological version. 

In the work \cite{PrepStatistic2019}, we develop a statistical framework on finding PCD based on observing average loss rates over a set of paths. We will discuss how imposing regularities based on sparsity and measurable PCD in the form of covariance among boundary set can be applied to deal with the highly underdetermined linear problem for the unknown link loss probabilities. An interesting future investigation would be to apply this form of PCD along with covariance based PCD to show the weighted based performance of proposed inference framework using network-level end-to-end measurements.

\section{Proofs of Theoretical Results}\label{sec:proofs}
\begin{proof}[\normalfont \textbf{Proof of Theorem~\ref{thm:ext}}]
  Introducing new notation
  \begin{equation*}
    \delta=\tcW - \cW
    \quad \mbox{and} \quad
    \eta = \cZ-A\cW,  
  \end{equation*}
  in \eqref{eq:matrix_const} and \eqref{eq:extpos:norm}, we seek
  solutions to the linear equation $A\delta = \eta$ of minimal
  $\ell_2$ norm $\Vert \delta \Vert^2$.  It is well known
  \cite{Penrose56}, that the least squares solution of minimal norm is given by $\delta = A'\eta$ where $A'$ is the Moore-Penrose inverse
  of $A$.
\end{proof}

\begin{proof}[\normalfont \textbf{Proof of Corollary~\ref{thm:inv:tree}}] 
  It suffices to show that $A$ has linearly independent rows.  Indeed, this implies that the null-space of $A^T$ is zero and, by the rank-nullity theorem, the range of $A$ is the entire space
  ($A \tcW = \cZ$ has a solution for any $\cZ$). Also, for any
  $\cW\neq 0$, $0 < \|A^T\cW\|^2 = \langle \cW, AA^T\cW \rangle$ and thus $AA^T \cW\neq 0$ and the matrix $AA^T$ is invertible.
  Since by our construction a path $P(b,b')$ does not pass through any other boundary vertices, every path is uniquely identified by the leaf it goes to (or from). Thereby every row of matrix $A$ has entry $1$ corresponding to the leaf edge. This entry must be $0$ in any other row in $A$. Therefore the rows are linearly independent.

\end{proof}

\begin{proof}[\normalfont \textbf{Proof of Corollary~\ref{cor:prune}}]
  Use Theorem~\ref{thm:ext} with $\cZ = A\cW$.  Full row rank of
  $\tilde A$ follows from Corollary~\ref{thm:inv:tree}.
\end{proof}

\begin{proof}[\normalfont \textbf{Proof of Lemma~\ref{lem:single_prune}}]	
	Denote the weight of $(u,v)$ by $\delta$.  Assume the first
	condition is satisfied (the proof of the second case is similar).
	We contract $(u,v)$, and adjust weights by subtracting $\delta$ from the weight of every edge originating from $u$ and adding $\delta$ to the weight of every edge terminating at $u$. It is easy to see that the weight of every path remains the same.
	Indeed, if the path $P$ visits the vertex $u$, $P$ must contain
	exactly one edge of the form $(u',u)$ and exactly one edge of the
	form $(u,v')$ (here we make use of the fact that there are no loops at vertex $u$). Then
	\begin{equation*}
	\widetilde W_{u',u} = W_{u',u} + \delta, \quad
	\widetilde W_{u,v'} = W_{u,v'} - \delta, \quad
	\widetilde W_{u',u} + \widetilde W_{u,v'} = W_{u',u} + W_{u,v'}.
	\end{equation*}
	The above reasoning also applies to the case $v'=v$ since we can
	view the contracted edge $(u,v)$ as an edge with $\widetilde W_{u.v}=0$. We have thus constructed one solution to
	\eqref{eq:constraint_after_prune} and therefore it is consistent.
\end{proof}

\begin{proof}[\normalfont \textbf{Proof of Theorem~\ref{thm:intrinsic}}]

	Writing $\tcW=\cW+\delta$, then we seek solutions $\delta$ to the
	linear equation $A\delta=-A\cW$ of minimal $L_2$ norm.  The form
	(\ref{eq:makec}) then follows if $AA^T$ is invertible, as we now
	establish. Observe we can write $A=(A_R,-A_S)$ (joining the
	corresponding rows) where $A_R$ and $A_S$ are the incidence matrices of the edges over paths in receiver trees and (reversed) paths in source trees.  Under an appropriate ordering of the pairs in $Q$ the matrix $A_R$ decomposes into sum $A_R=\oplus_{v\in V_b}A_{R,v}$ (similarly, $A_S=\oplus_{v\in V_b}A_{S,v}$ which we will need later). Corollary~\ref{thm:inv:tree} each $A_{R,v}$ has
	independent rows, therefore $A_R$ and, by extension, $A$ have
	independent rows.  Invertibility of $AA^T$ follows (see the proof of Corollary~\ref{thm:inv:tree} above).

	To prove the error estimate, equation~\eqref{eq:error_bound}, we
	write $\delta = \tcW-\cW=-A^T(AA^T)^{-1}A\cW$ and hence
	\begin{eqnarray*}
	\Vert\delta\Vert^2_2
	&=& \langle A^T(AA^T)^{-1}A\cW,A^T(AA^T)^{-1}A\cW\rangle\\
	&=& \langle (AA^T)(AA^T)^{-1}A\cW,(AA^T)^{-1}A\cW\rangle
	= \langle A\cW,(AA^T)^{-1}A\cW\rangle \\
	&\le& \Vert A\cW\Vert_2^2\, \Vert(AA^T)^{-1}\Vert    
	\end{eqnarray*}
	using the spectral norm for matrices. The result the follows if we can show $\Vert(AA^T)^{-1}\Vert\le 1/2$ or, equivalently,
	$\|AA^T\| \geq 2$.  The latter follows from
	$AA^T = A_RA_R^T + A_SA_S^T$ if
	$\langle y,A_RA_R^Ty\rangle \ge \Vert y\Vert^2$ and
	$\langle y,A_SA_S^Ty\rangle \ge \Vert y\Vert^2$ for all vectors
	$y\in\Rl^Q$. It suffices to show this for receiver trees. Observe
	$\langle y,A_RA_R^Ty\rangle = \sum_{v\in B} \langle y_v,
	A_{R,v}A_{R,v}^Ty_v\rangle$ where $y_v$ is projection of
	$y\in \Rl^Q$ onto $\Rl^{\cL_v}$.  We now represent
	\begin{equation*}
	A_{R,v}A_{R,v}^T = \sum_{i\in E^R_v} A_{R,v,i} A_{R,v,i}^T,
	\end{equation*}
	where $A_{R,v,i}$ is the $i$-th column of the matrix $A_{R,v}$.  If $i$ corresponds to an edge ending in a leaf $u$, the matrix
	$A_{R,v,i} A_{R,v,i}^T$ has a single 1 on the diagonal, in the
	position $((v,u),(v,u))$.  Therefore, the sum of
	$A_{R,v,i} A_{R,v,i}^T$ over leaf edges $i$ produces an identity
	matrix, while the rest of the summands are clearly positive
	semidefinite.  We conclude that
	$A_{R,v,i} A_{R,v,i}^T \geq \mathbb{I}$, therefore
	$\langle y_v, A_{R,v}A_{R,v}^Ty_v\rangle\ge \Vert y_v\Vert^2$ and so $\langle y,A_RA_R^Ty\rangle \ge \sum_{v\in V_B}\Vert y_v\Vert^2 = \Vert y\Vert^2$.

	We now show that the bound is tight for tree network graphs. Observe that, this class of networks has the property that the routing paths are symmetric and the receiver tree $\cT_v^R$ and source tree $\cT_v^S$ are isomorphic at each $v \in V_B$. 
	Denote by $S := A_SA_S^T$ and $R := A_RA_R^T$, these two matrices act on vectors indexed by pair of distinct boundary vertices. Due to structure of matrix $A$, $S_{(u_1,v_1)(u_2,v_2)} = 0$ if $u_1 \not = u_2$ and $S_{(u,v_1)(u,v_2)} = |P_{u,v_1} \cap P_{u,v_2}|$ i.e. number of directed edges in common between the two paths $P_{u,v_1}$ and $P_{u,v_2}$. If the routing is symmetric (this is the case on a tree network), then $R_{(v_1,u_1)(v_2,u_2)} = S_{(u_1,v_1)(u_2,v_2)}$ which imply $R = J^{-1} S J$ with $J$-permutation matrix sending $(u,v) \mapsto (v,u)$. Matrix $J$ satisfies $J_{(u_1,v_1)(u_2,v_2)} = \delta_{(u_1,v_2)}\delta_{(v_1,u_2)}$ with the property that $J^{-1} = J^T = J$. Since $AA^T = S + R$, therefore, it is sufficient to find $z$ such that 
	\begin{equation*}
	Sz = z, \quad \text{and} \quad SJz = Jz
	\end{equation*}    
	If $v_1$ and $v_2$ are siblings in $\cT_u^S$, then
	\begin{equation*}
	S_{(u,v_1)(u,v_1)} = S_{(u,v_2)(u,v_2)} = 1 + S_{(u,v_1)(u,v_2)} = 1 + S_{(u,v_2)(u,v_1)}
	\end{equation*} 
	Moreover, $\forall w \not = v_1, v_2$ 
	\begin{equation*}
	S_{(u,w)(u,v_1)} = S_{(u,w)(u,v_2)} 
	\end{equation*} 
	We conclude that $S-\mathbb{I}$ has identical columns $(u,v_1)$ and $(u,v_2)$ for any $v_1$ and $v_2$ siblings in $\cT_u^S$. Let now $u_1, u_2$ and $v_1, v_2$ be two distinct pairs of siblings (always possible to find these two set of distinct siblings in a tree with $|V_B| \geq 4$). Define vector $z$ with
	\begin{eqnarray*}
	z_{(u_1,v_1)} &=& +1, \quad z_{(u_2,v_1)} = -1 \\
	z_{(u_1,v_2)} &=& -1, \quad z_{(u_2,v_2)} = +1. 
	\end{eqnarray*}
	and all other $z_{(w_1,w_2)} = 0$. Since columns $(u_1,v_1)$ and $(u_1,v_2)$ of matrix $S-\mathbb{I}$ are identical and columns $(u_2,v_1)$ and $(u_2,v_2)$ are also identical, we get $(S-\mathbb{I})z = 0$, i.e. $Sz = z$.  
	
	But $\tilde z = Jz$ has the same structure
	\begin{eqnarray*}
	\tilde z_{(v_1,u_1)} &=& +1, \quad \tilde z_{(v_1,u_2)} = -1 \\
	\tilde z_{(v_2,u_1)} &=& -1, \quad \tilde z_{(v_2,u_2)} = +1. 
	\end{eqnarray*}
	Using properties of the columns $(v_1,u_1), (v_1,u_2)$ and $(v_2,u_1), (v_2,u_2)$ of matrix $S-\mathbb{I}$, we get $(S-\mathbb{I})\tilde z = 0$, i.e. $SJz = Jz$ as desired.  

	The above proof on tightness of bound holds for tree network graphs with $|V_B| \geq 4$. So, the only nontrivial tree which is not covered is 3-star one. Constructing incident matrix $A$ for the 3-star tree as discussed in Example \ref{exmp:intrinsicExample}, then $z = (+1,-1,-1,+1,+1,-1)^T$ is the eigenvector of matrix $AA^T$ with corresponding eigenvalue 2.  
	
\end{proof}


\section{Appendix: Application of Interior Point Methods for
  Positivity Constraints}\label{sec:appendix}
In this section we briefly discuss how interior-point method will be applied to find the minimizer of modified version of constrained quadratic optimization problem 
\begin{equation*}
  \label{eq:primalQPBarrier}
  x^* = \argmin_{x \in \mathbb{R}^n} \frac{1}{2} x^T H x + x^Tc - \mu \sum_{i=1}^{n} \ln (x_i) \quad \quad \text{subject to:} \quad  Bx = b,
\end{equation*}
appears in subsection \ref{sec:signConstraint} where $H \in \mathbb{R}^{n \times n}$ is positive semidefinite, 
$B \in \mathbb{R}^{m \times n}$ has full row rank, and parameter $\mu
\geq 0$ is log-barrier parameter (see \cite{AnLu07} for detailed discussion). Notice that in order to ensure that the objective function above is well defined, it is required that $x > 0$. The standard steps for solving this quadratic optimization problem as it will be discussed below, consists of (i) forming the Lagrangian, (ii) form the optimality condition, and finally (iii) apply iterative algorithm, for example, primal-dual path-following method. For the Lagrange multiplier
$\lambda \in \mathbb{R}^m$, the \textit{Lagrangian} is defined as
\begin{equation}
L(\alpha,\lambda,\mu) := \frac{1}{2} x^T H x+ c^T x - \mu \sum_{i=1}^{n} \ln (x_i) - \lambda^T (Bx - b) 
\end{equation}  
The conditions for a stationary point of Lagrangian with respect to $x$ and $\lambda$ satisfy 
\begin{equation}\label{eq:stationaryCond}
	\begin{split}
	Bx - b &= 0  \\
	B^T\lambda - \mu X^{-1} e -  H x - c &= 0 
	\end{split}
\end{equation}
for $x > 0$, where $X:= \text{diag}\{x_1, \dots, x_n \}$, and $e$ is vector of size $n$ with unit entries. If we let $s = \mu X^{-1} e$ to be a vector of size $n$, then $x > 0$ implies that $s > 0$ and equation \eqref{eq:stationaryCond} can be written as
\begin{equation}\label{eq:centralPath}
\begin{split}
Bx - b &= 0  \\
B^T\lambda + s - Hx - c &= 0 \\
Xs &= \mu e 
\end{split}
\end{equation}
for $x > 0$ and $s > 0$ which form the \textit{optimality condition}. Now in order to apply \textit{primal-dual path-following method}, let $w_k := \{x_k,\lambda_k,s_k\}$ be such that $x_k$ strictly feasible for the quadratic optimization problem ($x_k \in \mathbb{R}^n$ satisfies the constraints). The increment $\delta_w := \{\delta_x, \delta_{\lambda},\delta_s\}$ should be constructed such that the next iterate $w_{k+1} = w_k+\delta_w$ remains strictly feasible and approaches the central path. If $w_k$ were satisfy equation \eqref{eq:centralPath} with $\mu = \mu_{k+1}$, neglecting second-order term in equation \eqref{eq:centralPath}, we would have 
\begin{equation}\label{eq:centralPathItt}
\begin{split}
-H\delta_x + B^T\delta_{\lambda} + \delta_s &= 0 \\
B \delta_x &= 0 \\
S \delta_x + X \delta_s &= \mu_{k+1}e - X s_{k}
\end{split}
\end{equation}  
where $\Delta X := \text{diag}\{(\delta_x)_1 ,\ldots , (\delta_x)_n\}$ and $S := \text{diag}\{(s_k)_1 ,\ldots , (s_k)_n\}$. The solution to equation \eqref{eq:centralPathItt} then can be obtained as 
\begin{equation}\label{eq:centralPathSol}
\begin{split}
\delta_{\lambda} &= \Lambda y \\
\delta_x &= \Gamma X B^T \delta_{\lambda} - y \\
\delta_s &= H\delta_x - B^T \delta_{\lambda} 
\end{split}
\end{equation} 
where $\Gamma := (S + XH)^{-1}$, $\Lambda := (B \Gamma XB^T)^{-1}B$, and $y := \Gamma (X s_k - \mu_{k+1}e)$. In order to show that matrices $\Gamma$ and $\Lambda$ in equation \eqref{eq:centralPathSol} are well defined, observe that since $x_k > 0$ and $s_k > 0$, matrices $X$ and $S$ are positive definite. Therefore, the positive definite property of $X^{-1} S + H$ implies that the inverse of matrix $S + XH = X(X^{-1}S + H)$ exists. Applying the fact that $B$ is full row rank, then 
\begin{equation*}
	B\Gamma XB^T = B\big( X^{-1}S+H \big)^{-1} B^T
\end{equation*}
is also positive definite and hence nonsingular. As a result, this implies that $\delta_w$ has a unique solution. Once $\delta_w$ is calculated from equation \eqref{eq:centralPathSol}, for sufficiently small step size $\alpha_{k}$, the updated solution $w_{k+1} = w_k + \alpha_{k} \delta_w$ remains strictly feasible. In practice an appropriate choice of parameter $\mu_{k+1}$ is 
\begin{equation}
	\mu_{k+1} = \frac{x_k^T s_k}{n + \rho}
\end{equation}    
with $\rho \geq \sqrt{n}$. This would lead to an iteration complexity mentioned in subsection \ref{sec:signConstraint} for an $\varepsilon$-close solution to the optimal one.


\textbf{Acknowledgment}

We thank Kyriakos Manousakis for helpful discussions on the definition of evaluation metric in this paper.


\bibliographystyle{abbrv}
\bibliography{tomography} 

\end{document}